\documentclass[reqno]{amsart} 
\usepackage{fix-cm}
 \usepackage{fix-cm}
\usepackage{comment}
\usepackage{amsfonts,color}
\usepackage{latexsym,amssymb}
\usepackage{amsmath}
\usepackage[utf8]{inputenc}
\usepackage{enumitem}
\usepackage{xcolor}

\makeatletter
\renewcommand{\tagform@}[1]{\maketag@@@{\color{blue}(#1)}}
\makeatother

\newcommand{\Oe}{\Omega_{\ell}}
\newcommand{\Onl}{\Omega_{n\ell}}
\newcommand{\Lop}{\mathcal{L}}

\newcommand{\normaind}[2]{\left\lVert #1 \right\rVert_{#2}}
\newtheorem{theorem}{Theorem}[section]
\newtheorem{lemma}[theorem]{Lemma}
\newtheorem{proposition}[theorem]{Proposition}

\newtheorem{remark}[theorem]{Remark}
\newtheorem{definition}[theorem]{Definition}

\usepackage[bookmarksnumbered,colorlinks, linkcolor=blue, citecolor=red, pagebackref, bookmarks, breaklinks]{hyperref}
\usepackage{tikz-cd} 
\usepackage{xparse,tikz}
\usetikzlibrary{positioning,calc}
\usetikzlibrary {shapes.geometric}
\usetikzlibrary {patterns}
\usetikzlibrary {patterns,patterns.meta}
\usepackage[hyperpageref]{backref}

\catcode`@=11 \@addtoreset{equation}{section} \catcode`@=12
	
\begin{document}

\title[Coupled Local/Nonlocal Operators: Spectral Properties  
]{
Principal Eigenvalue for Coupled Local/Nonlocal Operators: Strict Monotonicity and Asymptotic Behavior }

\author[W. Cintra]{Willian Cintra$^{1}$}
\address{$^{1}$Universidade de Brasília, Departamento de Matem\'atica, BR-70910900 Brasília-DF, Brazil.  E-mails: willian@unb.br, jonatasperalta58@gmail.com}

	\author[J.  Peralta]{Jônatas Peralta$^1$}
    
    \author[A. Suárez]{Antonio Suárez$^{2}$}
\address{$^{2}$Departamento de Ecuaciones Diferenciales y Análisis Númérico, Universidad de Sevilla. E-mail: suarez@us.es}
		
\begin{abstract} 
In this paper, we investigate the properties of the generalized principal eigenvalue for a coupled local/nonlocal operator. First, we establish sufficient conditions to guarantee strict monotonicity properties, obtaining weaker hypotheses than those currently available in the literature for nonlocal operators. Next, we establish several results concerning the asymptotic behavior of this generalized eigenvalue, along with the convergence of the corresponding eigenfunctions. These results are very useful for the study of nonlinear problems. Our findings reveal novel phenomena that are not present in purely local or purely nonlocal cases. The main mathematical difficulty lies in the fact that the generalized eigenvalue does not always admit an eigenfunction, in addition to the inherent technicalities of the considered coupling.
\end{abstract}

\maketitle


\tableofcontents

\section{Introduction}
{Over the last decade, there has been a growing and intense interest in the mathematical study of reaction-diffusion models driven by mixed operators. These operators typically combine a classical local diffusion term, represented by the standard Laplacian ($-\Delta$), with a nonlocal dispersal operator governed by an integral kernel. Such models are particularly relevant in population dynamics and mathematical ecology. In this context, the movement of biological species is often not purely local; individuals exhibit a superposition of short-range random walks (modeled by the local Laplacian) and long-range dispersal strategies, such as seeds carried by wind or animal migrations (governed by the nonlocal integral kernel). 

In many scenarios, nonlocal diffusion is confined to a specific region of the habitat, whereas in the remainder of the domain, the phenomenon can be accurately governed by a purely local equation. To address such configurations, various coupled models have been proposed in the literature. See, for instance, \cite{delia2,MR4477810,MR4350298,DLLT,GW,Han,Kri} and references therein. In these works, the authors have explored the interplay between these two distinct diffusion regimes, revealing novel mathematical phenomena that neither operator exhibits independently. 

Recently, in \cite{RossiSIAM, ROSSICVPDE} (see also \cite{CAPANNAROSSIMJM,RosaJulio}), a rigorous mathematical coupling strategy was introduced, based on the minimization of an energy functional that consistently combines local and nonlocal interactions. Since then, such models have been intensively investigated from various analytical and applied perspectives, including evolution problems, homogenization, and numerical methods. See, for instance, \cite{AGBFR,delia3,MR4491890,PMLJ,RosaJulio}. 

More recently, in \cite{Peralta2026}, the authors established results concerning the existence of the principal eigenvalue and the maximum principle for this class of local/nonlocal coupled problems, which constitutes a fundamental step toward the study of positive solutions in nonlinear stationary models. 

Hence, the primary goal of this paper is to continue this line of research by exploring further properties of this eigenvalue problem and to rigorously characterize the asymptotic behavior of the principal eigenvalue and its corresponding eigenfunctions under determined perturbations. The limiting behavior turns out to be heavily dictated by the underlying geometry of the vanishing sets relative to the support of the nonlocal kernel.
}

\subsection*{Eigenvalue Problem and Hypotheses}
 
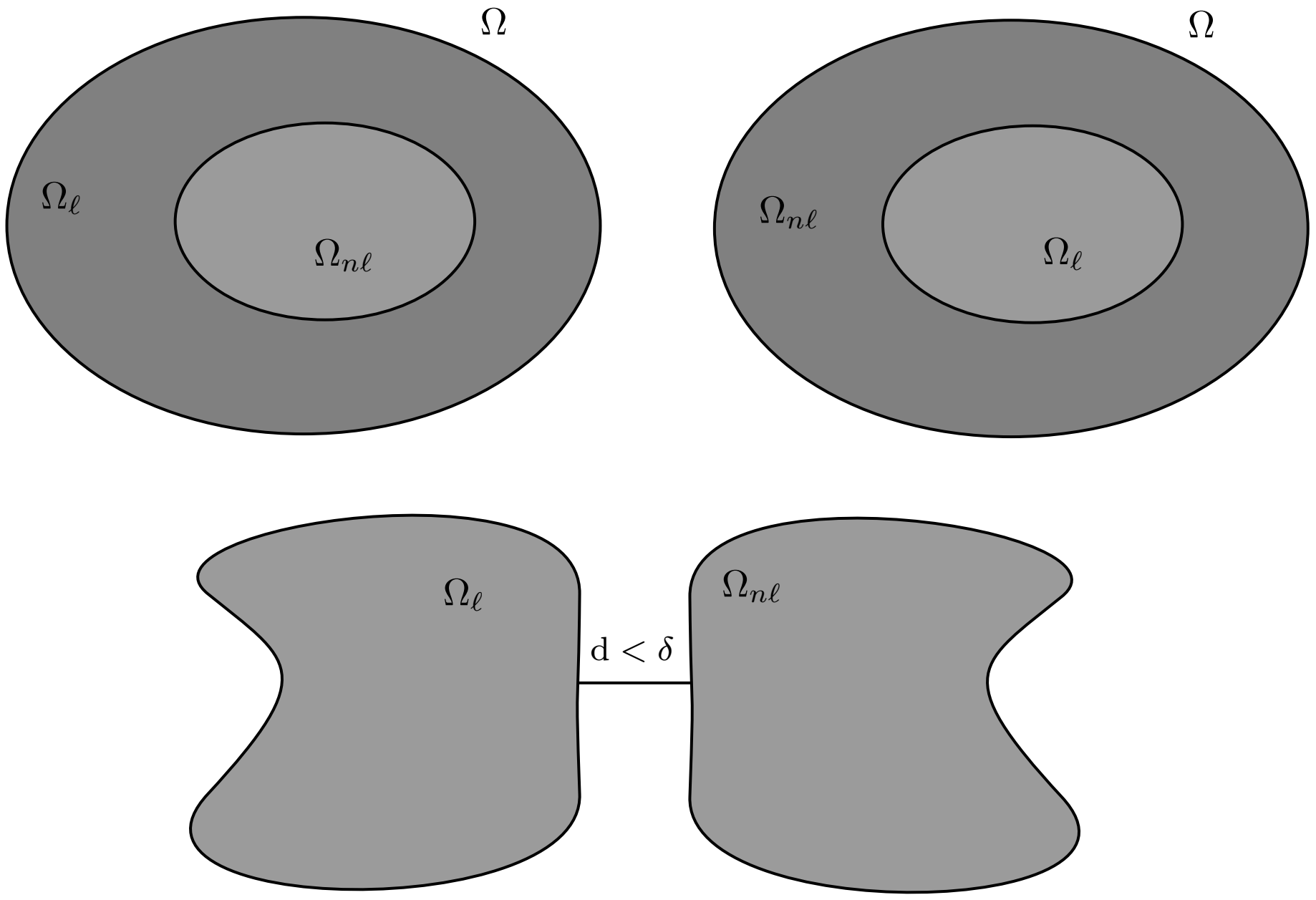
\begin{figure}
    \centering
    \tikzset{every picture/.style={line width=0.75pt}} 

\begin{tikzpicture}[x=0.75pt,y=0.75pt,yscale=-1,xscale=1]

\draw  [fill={rgb, 255:red, 128; green, 128; blue, 128 }  ,fill opacity=1 ] (43.67,99) .. controls (43.67,58.68) and (90.23,26) .. (147.67,26) .. controls (205.1,26) and (251.67,58.68) .. (251.67,99) .. controls (251.67,139.32) and (205.1,172) .. (147.67,172) .. controls (90.23,172) and (43.67,139.32) .. (43.67,99) -- cycle ;
\draw  [fill={rgb, 255:red, 155; green, 155; blue, 155 }  ,fill opacity=1 ] (102.67,97.5) .. controls (102.67,78.45) and (126.17,63) .. (155.17,63) .. controls (184.16,63) and (207.67,78.45) .. (207.67,97.5) .. controls (207.67,116.55) and (184.16,132) .. (155.17,132) .. controls (126.17,132) and (102.67,116.55) .. (102.67,97.5) -- cycle ;
\draw  [fill={rgb, 255:red, 128; green, 128; blue, 128 }  ,fill opacity=1 ] (291.67,100) .. controls (291.67,59.68) and (338.23,27) .. (395.67,27) .. controls (453.1,27) and (499.67,59.68) .. (499.67,100) .. controls (499.67,140.32) and (453.1,173) .. (395.67,173) .. controls (338.23,173) and (291.67,140.32) .. (291.67,100) -- cycle ;
\draw  [fill={rgb, 255:red, 155; green, 155; blue, 155 }  ,fill opacity=1 ] (350.67,98.5) .. controls (350.67,79.45) and (374.17,64) .. (403.17,64) .. controls (432.16,64) and (455.67,79.45) .. (455.67,98.5) .. controls (455.67,117.55) and (432.16,133) .. (403.17,133) .. controls (374.17,133) and (350.67,117.55) .. (350.67,98.5) -- cycle ;

\draw  [fill={rgb, 255:red, 155; green, 155; blue, 155 }  ,fill opacity=1 ] (244.42,228.02) .. controls (244.91,180.04) and (85.12,204.62) .. (114.09,228.02) .. controls (143.05,251.42) and (154.15,254.93) .. (114.09,298.22) .. controls (74.02,341.51) and (246,344.08) .. (244.42,298.22) .. controls (242.85,252.35) and (243.94,275.99) .. (244.42,228.02) -- cycle ;
\draw    (244,259.25) -- (284,259.25) ;
\draw  [fill={rgb, 255:red, 155; green, 155; blue, 155 }  ,fill opacity=1 ] (283.04,229.02) .. controls (282.56,181.04) and (442.34,205.62) .. (413.38,229.02) .. controls (384.42,252.42) and (373.31,255.93) .. (413.38,299.22) .. controls (453.45,342.51) and (281.47,345.08) .. (283.04,299.22) .. controls (284.62,253.35) and (283.53,276.99) .. (283.04,229.02) -- cycle ;

\draw (208,21.4) node [anchor=north west][inner sep=0.75pt]    {$\Omega $};
\draw (54,82.4) node [anchor=north west][inner sep=0.75pt]    {$\Omega _{\ell }$};
\draw (149.67,102.4) node [anchor=north west][inner sep=0.75pt]    {$\Omega _{n\ell }$};
\draw (456,22.4) node [anchor=north west][inner sep=0.75pt]    {$\Omega $};
\draw (405.17,101.9) node [anchor=north west][inner sep=0.75pt]    {$\Omega _{\ell }$};
\draw (305.67,87.4) node [anchor=north west][inner sep=0.75pt]    {$\Omega _{n\ell }$};
\draw (246.67,241.57) node [anchor=north west][inner sep=0.75pt]  [font=\fontsize{0.88em}{1.06em}\selectfont]  {$\operatorname{d< \delta }$};
\draw (195,221.4) node [anchor=north west][inner sep=0.75pt]    {$\Omega _{\ell }$};
\draw (292.67,218.4) node [anchor=north west][inner sep=0.75pt]    {$\Omega _{n\ell }$};

\end{tikzpicture}
    \caption{Admissible configurations of the domain. In the last case  $d:=\operatorname{dist}(\Oe,\Onl)<\delta$.}
    \label{fig: Dominio}
\end{figure}

Let us describe more precisely the problem under consideration. Let $\Omega_{\ell}$ and $\Omega_{n\ell}$ be disjoint regular subdomains of $\mathbb{R}^N$ ($N \geq 1$). We define the total domain as $\Omega = \Omega_{\ell} \cup \Omega_{n\ell} \cup \Gamma$, where $\Gamma := \partial \Omega_{\ell} \cap \partial \Omega_{n\ell}$ denotes their common boundary. Throughout this paper, we fix a constant $\delta > 0$ and assume that $$\operatorname{dist}(\Omega_{\ell}, \Omega_{n\ell}) < \delta.$$
Note that if the subdomains are strictly separated, $\Gamma$ is empty. However, this formulation naturally incorporates the scenario where $\Omega$ is partitioned by a hypersurface $\Gamma$; in this case, $\operatorname{dist}(\Omega_{\ell}, \Omega_{n\ell}) = 0$, which trivially satisfies the threshold condition (see Figure \ref{fig: Dominio}).

Since $\partial \Omega_\ell \neq \emptyset$ is composed of two disjoint components, $\Gamma_{0}:=\partial \Omega_{\ell} \cap \partial\Omega$ and $\Gamma=\partial \Omega_{\ell} \cap \Omega$, either of which may be empty, we consider the boundary operator $\mathfrak{B} : C(\Gamma_0) \times C^1(\Gamma) \to C(\partial \Omega_\ell)$ given by
\begin{equation*}
\mathfrak{B} u := \begin{cases}u & \text{on } \Gamma_0, \\ \frac{\partial u}{\partial \eta} & \text{on } \Gamma,\end{cases}
\end{equation*} where $\eta=\eta(x)$ denotes the outward unit normal field to $\Gamma=\partial\Oe \cap \Omega$.

We also introduce the following Banach spaces:
\begin{align*}
    X&= \{u \in L^2(\Omega) : u|_{\Omega_\ell} \in W^{2,p}(\Omega_\ell),\, p > N, \text{ and }  u|_{\overline{\Omega}_{n\ell}} \in C(\overline{\Omega}_{n\ell}) \};   \\ X_{\mathfrak{B}}&= \{u \in X : \, \mathfrak{B}u=0 \text{ on } \partial \Omega_{\ell} \};\\
    W_\mathfrak{B}^{2,p}(\Omega_{\ell}) &= \{ u \in W^{2,p}(\Omega_{\ell}): \mathfrak{B}u=0\} \quad \mbox{and}\\
    Y&=\{u \in L^2(\Omega): \, u|_{\Oe} \in L^p(\Oe),\, p>N, \text{ and } u|_{\Onl} \in C(\overline{\Omega}_{n\ell})\}.   \end{align*}
Motivated by the models introduced in \cite{RossiSIAM, ROSSICVPDE} and the recent results obtained in \cite{Peralta2026}, in this paper we establish new monotonicity properties and analyze the asymptotic behavior of the principal eigenvalue of the following eigenvalue problem:

\[
\begin{cases} \tag{$P$}
\Lop_1 u + v_1(x)u= \lambda u  & \text{in } \Omega_{\ell}, \\ 
\mathfrak{B}u = 0 & \text{on } \partial \Omega_{\ell}, \\ 
\Lop_2 u  +v_2(x)u=\lambda u & \text{in } \Omega_{n \ell}, \label{P} 
\end{cases} 
\]
where $v_{1}\in L^\infty(\Oe)$, $v_{2} \in C(\overline{\Omega}_{n\ell})$  and $\mathcal{L}_1: X \to L^2(\Omega_\ell)$, $\mathcal{L}_2: X \to L^2(\Omega_{n\ell})$ are given by
\begin{align*}
\mathcal{L}_1u(x) &= -\Delta u - \int_{\Omega_{n\ell}} J(x-y)u(y)dy \\
\mathcal{L}_2u(x) &= -2\int_{\Omega_{n\ell}} J(x-y)u(y)dy - \int_{\Oe} J(x-y)u(y)dy.
\end{align*}
For the dispersion kernel $J$ we assume that 
\begin{enumerate}[label=$(J_{\arabic*})$]\item \label{J1} $J \in C(\mathbb{R}^N)$ is a nonnegative function such that $J(0)>0$ and there exists $\delta>0$ such that $\operatorname{supp}J = B_\delta(0)$.
\end{enumerate}
Note that we can make the natural identification of $X_\mathfrak{B}$ with $W_\mathfrak{B}^{2,p}(\Omega_\ell) \times C(\overline{\Omega}_{n\ell})$ such that each $w \in X_\mathfrak{B}$ is identified with the pair $(u,v)\in W_\mathfrak{B}^{2,p}(\Omega_\ell) \times C(\overline{\Omega}_{n\ell})$ where $u=w\vert{}_{\Oe}$ and $v=w\vert{}_{\Onl}$. Analogously, $Y$ is identified with $L^p(\Oe)\times C(\overline{\Omega}_{n\ell})$. Thus, defining $V= (v_1,v_2)\in L^\infty(\Oe)\times C(\overline{\Omega}_{n\ell})$ and $\mathfrak L:X\to Y$, where
\[
 \mathfrak{L}u:= \begin{cases} \Lop_1 u  & \text{in } \Omega_{\ell}, \\ \Lop_2 u  & \text{in } \Omega_{n \ell},\end{cases}
\]
the eigenvalue problem \eqref{P} can be rewritten as
\[
\begin{cases} 
\mathfrak{L}u +V(x)u = \lambda u & \text{in } \Omega,\\ 
\mathfrak{B} u=  0& \text{on } \partial \Omega_{\ell}.

\end{cases}
\] 
Inspired by the works of \cite{BNV,MR2718672,DV75}, we introduce the set
$$ \mathfrak{P}= \{ u \in X:  u(x)>0 \mbox{ for all } x \in \Omega \mbox{ and } \mathfrak{B}u \geq 0 \mbox{ on } \partial\Oe\}$$
and define the so-called generalized eigenvalue and the constant $m$ as
$$ \lambda_p^{\Omega}[\mathfrak{L}+V] := \sup_{\phi \in \mathfrak{P}} \inf_{x \in \Omega} \frac{(\mathfrak{L} + V(x))\phi(x)}{\phi(x)} \quad \mbox{and} \quad m:= \min_{x \in \overline{\Omega}_{nl}}  v_{2}(x).$$
It was shown in \cite{Peralta2026} that $\lambda_{p}^\Omega[\mathfrak{L}+V]$ is a principal eigenvalue of \eqref{P} if and only if 
$$ \lambda_p^{\Omega}[\mathfrak{L}+V]< m. $$ 
Furthermore, several other important properties were established (see also Section~\ref{sec2.3}).

\subsection*{Main Results} This work has two main goals. The first is to establish sufficient conditions for the strict monotonicity of $\lambda_p^\Omega[\mathfrak{L}+V]$ with respect to $V$ and the domain $\Omega$, as well as to relate this principal eigenvalue to the principal eigenvalues of the local and nonlocal domains (see Section \ref{sec3}).

It is worth pointing out that strict monotonicity always holds for the principal eigenvalue of uniformly elliptic (local) operators, whereas this is not necessarily the case for integral operators. To the best of our knowledge, only a few partial results addressing this issue are currently available in the literature. For instance, in \cite{MR3485862}, the authors require both principal eigenvalues to have a positive eigenfunction. We also mention \cite{WANG2027110124}, which considers the one-dimensional case.

Our findings show that in both coupled and nonlocal settings, it suffices for one of the generalized eigenvalues to admit an eigenfunction. This result is very useful for establishing non-existence results in nonlinear problems.

Next, we analyze several asymptotic behaviors of this generalized eigenvalue.  It is worth noting that such limiting properties are well known in the purely local setting (see \cite{DuH,FRAILE1996,HessKato,LG-JDE96}) as well as for certain integral operators (see \cite{WANG2027110124}), and it allows for the analysis of the structure of positive solutions in logistic models with refuges. However, due to the coupling under consideration, various configuration possibilities arise, revealing novel phenomena. Furthermore, our proofs are entirely original and differ fundamentally from the approaches presented in the aforementioned works.

To state our results precisely, let $A := (a_1, a_2) \in L^{\infty}(\Oe) \times C(\overline{\Omega}_{n\ell})$, $\gamma_{k}:=(\alpha_{k}, \beta_{k}) \in C(\overline{\Omega}_{\ell} )\times C(\overline{\Omega}_{n\ell})$, $k \ge 1$, and define $V:= A + \gamma_k$. Adopting the notation 
$$\mu_{k}:= \lambda_{p}^{\Omega}[\mathfrak{L}+A+\gamma_{k}],$$
we have the following result, which analyzes the behavior as $\gamma_k \to -\infty$. 

\begin{proposition}\label{prop-gammak} Suppose that there exists an open set $D$ such that either
     \begin{enumerate}
         \item [\rm (i)] $D \subset \Oe$  and
            \begin{equation} \label{alphauniK}
                \sup_D\alpha_{k}(x) \to -\infty,
            \end{equation}
    or
    \item[\rm (ii)]  $D \subset \Onl$  and 
\begin{equation} \label{betauniK}
    \sup_D\beta_{k}(x) \to -\infty.
\end{equation}
Then
    $$\mu_{k} \to  -\infty.$$
     \end{enumerate}
 \end{proposition}
To state the next result, we need to introduce some additional notation. Let $D \subset \mathbb{R}^N$ be a bounded domain and $a \in C(\overline{D})$. Setting
$$L_Du:= -\int_D J(x-y) u(y) \,dy \quad \forall u \in C(\overline{D}).$$we define
\begin{equation*}
    \lambda_p^D[L_D+a] := \sup_{u\in C(\overline{D}), u>0} \inf_{x \in D} \frac{(L_D+a(x))u(x)}{u(x)}.
\end{equation*}
Moreover, when $D$ is a regular bounded domain, for each $a \in L^\infty(D)$ we also denote by $\sigma_1^D[-\Delta+a]$ the principal eigenvalue of  $-\Delta+a$ in $H_0^1(D)$. In particular, defining
\[
L_{n\ell}u:=-2 \int_{\Onl}J(x-y) u(y) dy \quad \forall u \in C(\overline{\Omega}_{n\ell}),\]
we will consider $\lambda_{p}^{\Onl}[L_{n\ell}+a_{2}]$. These eigenvalues are well known in the literature; see, for instance, \cite{MR3583503} and \cite{LG-JDE96}. See also Section \ref{Secao2}.

Thus, we have the following result:

\begin{theorem}\label{T1}
   Let $\gamma_k:=(\alpha_{k},\beta_{k}) \in C(\overline{\Omega}_{\ell}) \times C(\overline{\Omega}_{n\ell})$ be a non-decreasing sequence and $(\alpha^*,\beta^*) \in C(\overline{\Omega}_\ell)\times C(\overline{\Omega}_{n\ell})$.
    \begin{enumerate}
        \item[(i)] If $\alpha_k(x) \to +\infty$ for all $x \in \overline{\Omega}_\ell$ and $\beta_k \to \beta^*$ in $C(\overline{\Omega}_{n\ell})$, then
        $$ \mu_{k} \to\lambda_{p}^{\Onl}[L_{n\ell}+a_2+\beta^*].$$
        \item[(ii)] If $\alpha_k \to \alpha^*$ in $C(\overline{\Omega}_\ell)$ and $\beta_k(x) \to +\infty$ for all $x \in \overline{\Omega}_{n\ell}$, then
    $$\mu_{k} \to \sigma_{1}^{\Oe}[-\Delta+a_1+\alpha^*;\mathfrak{B}].$$
    \item[(iii)] If $\alpha_k(x) \to +\infty$ for all $x \in \overline{\Omega}_\ell$ and $\beta_k(x) \to +\infty$ for all $x \in \overline{\Omega}_{n\ell}$, then
    $$\mu_{k} \to +\infty.$$   \end{enumerate}
\end{theorem}

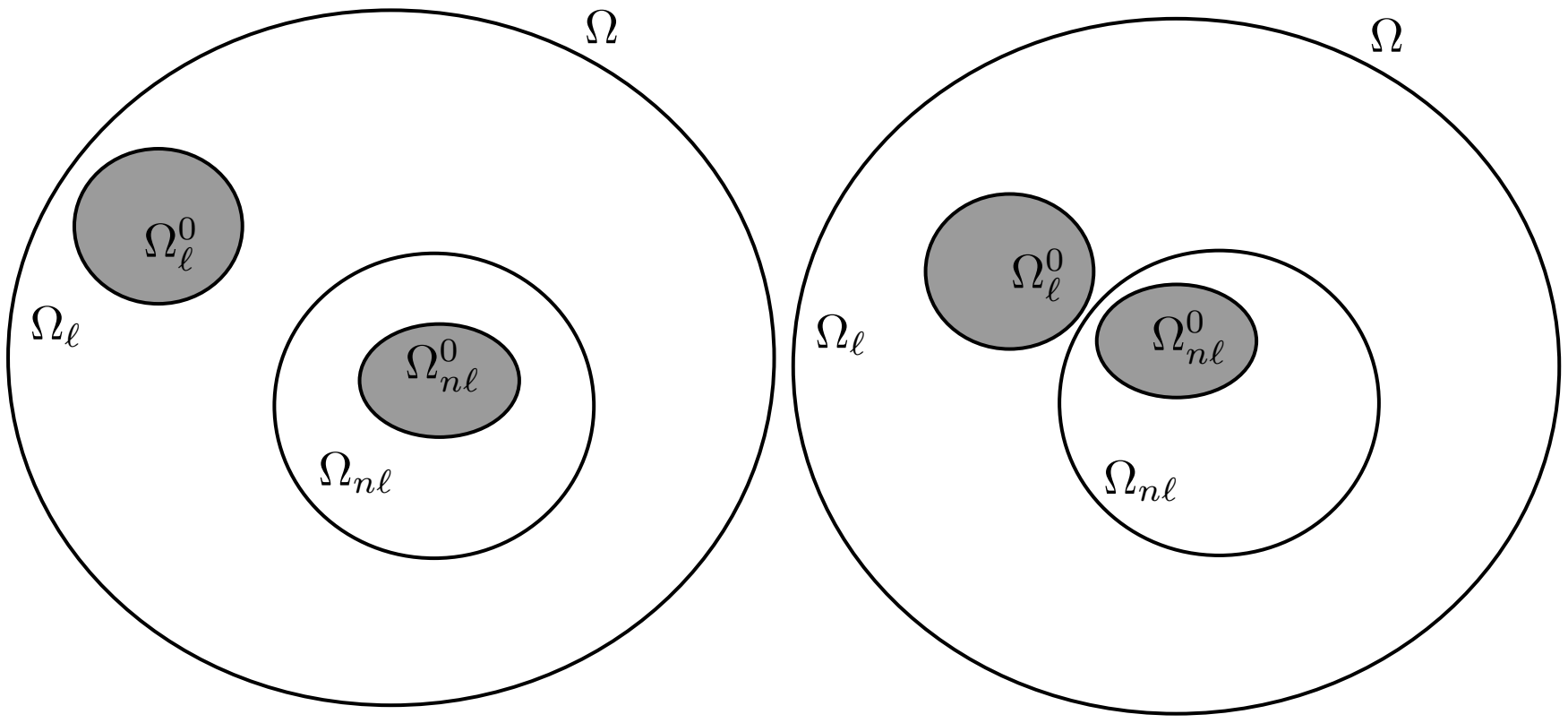
\begin{figure}
    \tikzset{every picture/.style={line width=0.75pt}} 

\begin{tikzpicture}[x=0.75pt,y=0.75pt,yscale=-1,xscale=1]

\draw  [fill={rgb, 255:red, 255; green, 255; blue, 255 }  ,fill opacity=1 ] (4,190.12) .. controls (4,136.36) and (52.06,92.78) .. (111.33,92.78) .. controls (170.61,92.78) and (218.67,136.36) .. (218.67,190.12) .. controls (218.67,243.88) and (170.61,287.46) .. (111.33,287.46) .. controls (52.06,287.46) and (4,243.88) .. (4,190.12) -- cycle ;
\draw  [fill={rgb, 255:red, 255; green, 255; blue, 255 }  ,fill opacity=1 ] (78.61,203.58) .. controls (78.61,180) and (98.65,160.88) .. (123.38,160.88) .. controls (148.1,160.88) and (168.15,180) .. (168.15,203.58) .. controls (168.15,227.16) and (148.1,246.28) .. (123.38,246.28) .. controls (98.65,246.28) and (78.61,227.16) .. (78.61,203.58) -- cycle ;
\draw  [fill={rgb, 255:red, 155; green, 155; blue, 155 }  ,fill opacity=1 ] (102.49,196.52) .. controls (102.49,187.77) and (112.51,180.68) .. (124.87,180.68) .. controls (137.23,180.68) and (147.25,187.77) .. (147.25,196.52) .. controls (147.25,205.27) and (137.23,212.36) .. (124.87,212.36) .. controls (112.51,212.36) and (102.49,205.27) .. (102.49,196.52) -- cycle ;
\draw  [fill={rgb, 255:red, 155; green, 155; blue, 155 }  ,fill opacity=1 ] (22.55,153.29) .. controls (22.55,141.3) and (33.09,131.58) .. (46.1,131.58) .. controls (59.11,131.58) and (69.66,141.3) .. (69.66,153.29) .. controls (69.66,165.28) and (59.11,175.01) .. (46.1,175.01) .. controls (33.09,175.01) and (22.55,165.28) .. (22.55,153.29) -- cycle ;
\draw  [fill={rgb, 255:red, 255; green, 255; blue, 255 }  ,fill opacity=1 ] (224,192.49) .. controls (224,138.73) and (272.05,95.15) .. (331.33,95.15) .. controls (390.61,95.15) and (438.67,138.73) .. (438.67,192.49) .. controls (438.67,246.25) and (390.61,289.83) .. (331.33,289.83) .. controls (272.05,289.83) and (224,246.25) .. (224,192.49) -- cycle ;
\draw  [fill={rgb, 255:red, 255; green, 255; blue, 255 }  ,fill opacity=1 ] (298.61,202.79) .. controls (298.61,179.21) and (318.65,160.09) .. (343.38,160.09) .. controls (368.1,160.09) and (388.14,179.21) .. (388.14,202.79) .. controls (388.14,226.37) and (368.1,245.49) .. (343.38,245.49) .. controls (318.65,245.49) and (298.61,226.37) .. (298.61,202.79) -- cycle ;
\draw  [fill={rgb, 255:red, 155; green, 155; blue, 155 }  ,fill opacity=1 ] (309.06,185.43) .. controls (309.06,176.69) and (319.08,169.59) .. (331.44,169.59) .. controls (343.8,169.59) and (353.82,176.69) .. (353.82,185.43) .. controls (353.82,194.18) and (343.8,201.27) .. (331.44,201.27) .. controls (319.08,201.27) and (309.06,194.18) .. (309.06,185.43) -- cycle ;
\draw  [fill={rgb, 255:red, 155; green, 155; blue, 155 }  ,fill opacity=1 ] (261.09,165.96) .. controls (261.09,153.97) and (271.64,144.25) .. (284.65,144.25) .. controls (297.66,144.25) and (308.2,153.97) .. (308.2,165.96) .. controls (308.2,177.96) and (297.66,187.68) .. (284.65,187.68) .. controls (271.64,187.68) and (261.09,177.96) .. (261.09,165.96) -- cycle ;

\draw (164.02,91.41) node [anchor=north west][inner sep=0.75pt]    {$\Omega $};
\draw (8.63,174.35) node [anchor=north west][inner sep=0.75pt]    {$\Omega _{\ell }$};
\draw (89.41,215.22) node [anchor=north west][inner sep=0.75pt]    {$\Omega _{n\ell }$};
\draw (113.71,183.54) node [anchor=north west][inner sep=0.75pt]    {$\Omega _{n\ell }^{0}$};
\draw (40.42,149.49) node [anchor=north west][inner sep=0.75pt]    {$\Omega _{\ell }^{0}$};
\draw (384.01,93.79) node [anchor=north west][inner sep=0.75pt]    {$\Omega $};
\draw (228.63,176.73) node [anchor=north west][inner sep=0.75pt]    {$\Omega _{\ell }$};
\draw (309.4,217.6) node [anchor=north west][inner sep=0.75pt]    {$\Omega _{n\ell }$};
\draw (322.83,175.62) node [anchor=north west][inner sep=0.75pt]    {$\Omega _{n\ell }^{0}$};
\draw (283.45,158.2) node [anchor=north west][inner sep=0.75pt]    {$\Omega _{\ell }^{0}$};

\end{tikzpicture}
    
     \caption{Possible configuration of subdomains }\label{fig:2}
\end{figure}

Next, let us analyze the case where the weight functions may be identically zero in certain subdomains. To this end, consider $B=(b_1,b_2) \in C(\overline{\Omega}_\ell)\times C(\overline{\Omega}_{n\ell})$ such that $b_1,b_2\geq0$ and define:
$$
\Omega_{\ell}^0= \operatorname{int}\{ x \in \Oe; \,b_1(x)=0\} \quad \mbox{and} \quad
\Omega_{n\ell}^0= \operatorname{int}\{ x \in \Onl; \,b_2(x)=0\}.
$$
Let us assume that these sets are regular subdomains, which may be empty, and $\overline{\Omega}_\ell^0 \subset \Oe$. When both are non-empty, we define 
$$\Omega_0=\Omega_{\ell}^0 \cup \Omega_{n\ell}^0 \subset \Omega.$$
Figure \ref{fig:2} illustrates two admissible configurations for $\Oe^{0}$ and $\Onl^{0}.$ 

Thus, we can define a local/nonlocal operator on $\Omega_0$ as follows. Let
$$X_0= \{u \in L^2(\Omega_0) : u\vert{}_{\Omega_{\ell}^0} \in W_0^{2,p}(\Omega_{\ell}^0),\, p > N, \text{ and }  u\vert{}_{\overline{\Omega}_{n\ell}^0} \in C(\overline{\Omega}_{n\ell}^0) \},$$
$$Y_0=\{u \in L^2(\Omega_0): \, u\vert{}_{\Omega_{\ell}^0} \in L^p(\Oe^0),\, p>N, \text{ and } u\vert{}_{\overline{\Omega}_{n\ell}^0} \in C(\overline{\Omega}_{n\ell}^0)\}.$$
We introduce the mapping $\mathfrak{L}_{0}:X_0 \to Y_0$ given by 
$$\mathfrak{L}_{0}u:= \begin{cases} \mathcal{L}_1^{0} u  & \text{in } \Omega_{\ell}^{0}, \\ \mathcal{L}_2^{0} u  & \text{in } \Omega_{n \ell}^{0}, \end{cases}$$
where
$$\begin{aligned} \mathcal{L}_1^{0}u(x)& := -\Delta u - \int_{\Omega_{n\ell}^{0}} J(x-y)u(y)dy, \\ \mathcal{L}_2^{0}u(x) &:= -2\int_{\Omega_{n\ell}^0} J(x-y)u(y)dy - \int_{\Oe^0} J(x-y)u(y)dy. \end{aligned}$$
Finally, let us also consider the following operator $L_{n\ell}^{0}:C(\overline{\Omega}_{n\ell}^{0}) \to C(\overline{\Omega}_{n\ell}^{0})$ given by $$L_{n\ell}^{0}u :=-2\int_{\Omega_{n\ell}^{0}}J(x-y) u(y)\, dy.$$
Thus, if for each constant $\gamma>0$ we define 
$$
\Sigma(\gamma):= \lambda_p^\Omega[\mathfrak{L}+A+\gamma B],
$$
then the following result holds.

\begin{theorem}\label{T2}
Suppose that  $\Oe^0 \neq \emptyset$ and $\Onl^{0} \neq \emptyset$. Then
\begin{equation}
\lim_{\gamma\to+\infty} \Sigma(\gamma)= \begin{cases} 
\lambda_p^{\Omega_0}[\mathfrak{L}_0+A] & \text{if } \operatorname{dist}(\Oe^0,\Onl^0)<\delta, \\
\min \left\{\sigma_{1}^{\Omega_{\ell}^0}[-\Delta +a_1], \lambda_p^{\Omega_{n\ell}^0}[L_{n\ell}^{0}+a_2]\right\} & \text{if } \operatorname{dist}(\Oe^0,\Onl^0)\geq\delta.
\end{cases}
\end{equation}
\end{theorem}

\begin{remark}
\begin{enumerate}   
\item Theorem \ref{T2} yields a similar result when $\Oe^0$ or $\Onl^{0}$ (or both) are empty (see Section \ref{sec4.2}). Specifically, it holds that: 
    \begin{equation}
        \lim_{\gamma\to+\infty} \Sigma(\gamma)= \begin{cases} 
\sigma_{1}^{\Oe^{0}}[-\Delta+a_{1}] & \text{if } \Omega_{\ell}^{0}\neq\emptyset \text{ and } \overline{\Omega}_{n\ell}^{0}=\emptyset, \\
\lambda_{p}^{\Onl^{0}}[L_{n\ell}^{0}+a_{2}] & \text{if } \overline{\Omega}_{\ell}^{0}=\emptyset \text{ and } \Omega_{n\ell}^{0}\neq\emptyset, \\
+\infty & \text{if } \overline{\Omega}_{\ell}^{0}=\emptyset \text{ and } \overline{\Omega}_{n\ell}^{0}=\emptyset.
\end{cases}
\end{equation}
\item 
It is worth noting that we also prove the convergence of the corresponding eigenfunctions in Theorems \ref{T1} and \ref{T2}, provided that the limit values possess an associated positive eigenfunction. See Section \ref{sec4}. 
\item   With minor adaptations, our arguments can be used to prove Theorem \ref{T2} when $\overline{\Omega}_\ell^0 \subset \Oe$ is replaced by $\Oe^0 \subset \Oe$ and $\Gamma \subset \partial \Oe^0$.

\end{enumerate}
\end{remark}

We emphasize the importance of Theorem \ref{T2} and its applicability, particularly in the analysis of the asymptotic behavior of degenerate logistic models. This type of argument appears in the literature, for instance, in the context of local \cite{FRAILE1996,SJW}, nonlocal \cite{MR3485862}, and periodic parabolic problems \cite{Anjos2026,MR4867294}. Furthermore, this theorem will allow us to establish novel behaviors in population dynamics models, which will be explored in detail in a forthcoming paper.



Finally, for the next result, we let $\delta \to 0$ to study the asymptotic behavior of $\lambda_p^\Omega[\mathfrak{L}+V]$.  Since hypothesis \ref{J1} states that $\operatorname{supp}J=B_{\delta}(0)$, it is convenient to denote
$$J_{\delta}=J, \quad \mathfrak{L}_{\delta}=\mathfrak{L} \quad \text{and } \lambda(\delta):=\lambda_{p}^{\Omega}[\mathfrak{L}_{\delta}+A].$$
Thus, we obtain the following result:

\begin{theorem}\label{T3}
    
 Assume that there exist $\delta', C > 0$ such that 
   \begin{equation} \label{J-limitada}
       |J_{\delta}(z)| \le C \quad \forall \delta<\delta', \, z \in \mathbb{R}^N.
   \end{equation} 
   If $\operatorname{dist}(\Oe,\Onl)=0$, then 
   \begin{equation}\label{delta-to-0}
       \lambda(\delta) \to  \min \left\{\sigma_{1}^{\Oe}[-\Delta+a_{1};\mathfrak{B}], \min_{\overline{\Omega}_{n\ell}} a_2(x)\right\}\quad \text{as } \delta \to 0^+.
   \end{equation}

\end{theorem}
\subsection*{Additional Notations} Let us 
introduce some notations that will be used throughout the paper.
\begin{itemize}\item $P_\ell:=P_{C_{0}^{1}(\overline{\Omega}_{\ell})}$ and $P_{n\ell}:=P_{C(\overline{\Omega}_{n\ell})}$ are the positive cones of $C_{0}^{1}(\overline{\Omega}_{\ell})$ and $C(\overline{\Omega}_{n\ell})$ respectively.
\item When $u \in X$, we write $u\gg 0$  to mean that
$$u\vert{}_{\Omega_\ell} \in \operatorname{int} (P_{\ell}) \quad \mbox{and} \quad u\vert{}_{\overline{\Omega}_{n\ell}} \in \operatorname{int} (P_{n\ell}).$$
\item When $u \in Y$, we write $u>0$  to mean that $u \geq 0$ a.e. in $\Omega$ and $u \neq 0$. In this case we say that $u$ is {\it positive}. Similarly, we write $u>0$ in $\Onl$ or $u>0$ in $\Oe$.
\item $\vert{}\cdot\vert{}_{2}$ denotes the usual norm in $L^{2}(\Omega)$.
\item   $H^1_{\Gamma_0}(\Oe) = \{ u \in H^1(\Oe) :  u = 0 \text{ on } \Gamma_0 \}.$
\item Let \(L_{n\ell}:C(\overline{\Omega}_{n\ell})\to C(\overline{\Omega}_{n\ell})\) and \(L_\ell:X\to C(\overline{\Omega}_{n\ell})\) be defined by
\[
L_{n\ell}u:=-2\int_{\Omega_{n\ell}}J(x-y)u(y)\,dy,
\qquad
L_\ell u:=-\int_{\Omega_\ell}J(x-y)u(y)\,dy.
\]
\end{itemize} 

This paper is organized as follows. In Section \ref{Secao2}, we review auxiliary results from the literature regarding local and nonlocal eigenvalues, as well as the eigenvalue $\lambda_{p}^{\Omega}[\mathfrak{L}+V]$. Section \ref{sec3} establishes sufficient conditions to obtain strict monotonicity with respect to $V$ and the domain $\Omega$, also detailing the relationship of $\lambda_{p}^{\Omega}[\mathfrak{L}+V]$ with the local and nonlocal eigenvalues. Section \ref{sec4} is dedicated to the study of the asymptotic behavior of the eigenvalue, where we will prove Proposition \ref{prop-gammak} and Theorems \ref{T1} and \ref{T2}. Finally, the analysis of the behavior with respect to $\delta$ is carried out in Section \ref{sec5}.

\section{Preliminaries}\label{Secao2}

This section is dedicated to reviewing key properties of three eigenvalue problems that are pivotal to the present work. For the reader's convenience, we divide this section into three parts. Sections \ref{sec2.1} and \ref{sec2.2} revisit a nonlocal and a local eigenvalue problem, respectively. Section \ref{sec2.3} introduces a local/nonlocal eigenvalue problem associated with equation \eqref{P}. 
\subsection{Nonlocal eigenvalue problem} \label{sec2.1}
Let $D \subset \mathbb{R}^N$ be a bounded  domain and $a \in C(\overline{D})$. Setting
$$L_Du:= -\int_D J(x-y) u(y) \,dy,$$
we define
\begin{equation}\label{autovalor nonlocal}
    \lambda_p^D[L_D+a] := \sup_{u\in C(\overline{D}), u>0} \inf_{x \in D} \frac{(L_D+a(x))u(x)}{u(x)}.
\end{equation}
Thus, we have \begin{equation}\label{lambdap-NL<m}\lambda_p^D[L_D+a]\leq \min_{\overline{D}}a.\end{equation}The following result collects some properties of this generalized eigenvalue that will be used throughout the paper.

\begin{lemma}\label{propriedades-NL-autovalor}The following statements hold.
\begin{enumerate}\item[\rm(i)] $\lambda_p^D[L_D+a]$ is continuous and nondecreasing with respect to $a\in C(\overline{D})$.
\item[\rm(ii)] If $D_1 \subset D$ is a regular domain, then
    
    $$\lambda_p^{D}[L_{D}+a]\leq \lambda_p^{D_1}[L_{D_1}+a]$$
    \item[\rm(iii)] $\lambda_p^D[L_D+a]$ is an eigenvalue of $L_D+a$ with positive eigenfunction $\varphi_D \in C(\overline{D})$ if and only if $$\lambda_p^D[L_D+a] < \min_{\overline{D}}a.$$
\end{enumerate}
\end{lemma}\begin{proof} See, for instance, \cite{BCV,MR2718672,MR3583503,MR3485862}.\end{proof}

\begin{remark}
To the best of our knowledge, the fact that the monotonicity provided by items (i) and (ii) is strict is known only in a few particular cases (see \cite{WANG2027110124} and \cite{MR3485862}). In Section \ref{sec3} we will prove that for the coupled local/nonlocal operator, it is sufficient that only one of them admits an eigenfunction. Moreover, our arguments can be used to establish the same result for these nonlocal operators.
\end{remark}

The following result will also be useful throughout the paper.
\begin{lemma}
\label{solução trivial local}
Suppose that  there exists a function $\phi \in C(\overline{D})$ such that $\phi(x)>0$ for all $x \in \overline{D}$ and
$$
(L_D+a)\phi >0\quad  \mbox{in }\overline{D}.
$$
Then the homogeneous linear problem
\begin{equation}
   L_{D}u+au=0 \quad \mbox{in }  \overline{D}
\end{equation}
possesses only the trivial solution $u\equiv0$ in $C(\overline{D})$. 
\end{lemma}
\begin{proof}
    It is a consequence of \cite[Lemma 2.1]{MR3583503} by using the same arguments as in \cite[Corollary 2.5]{Peralta2026}.
\end{proof}

\begin{remark} 
Recall that the operator $L_{n\ell}:C(\overline{\Omega}_{n\ell}) \to C(\overline{\Omega}_{n\ell})$ is defined by
$$
L_{n\ell}u(x):= -2\int_{\Onl}J(x-y) u(y)\, dy.
$$
Thus, all the properties stated in Lemma \ref{propriedades-NL-autovalor} and Lemma \ref{solução trivial local} remain valid for $L_{n\ell}+a_2$ and its corresponding generalized principal eigenvalue $\lambda_{p}^{\Onl}[L_{n\ell}+a_{2}]$.
\end{remark}

\subsection{Local eigenvalue problem}\label{sec2.2}
Let $c\in L^\infty(D)$, where $D \subset \mathbb{R}^N$ is a regular subdomain such that $\partial D= \Gamma_{0,D} \cup \Gamma_{1,D}$, with $\Gamma_{0,D}$ and $\Gamma_{1,D}$ being disjoint open and closed subsets of $\partial D$. We define\begin{equation*}\mathfrak{B}_1 u := \begin{cases}
u & \text{on } \Gamma_{0,D}, \\ \frac{\partial u}{\partial \eta} & \text{on } \Gamma_{1,D},\end{cases}\end{equation*}where $\eta=\eta(x)$ denotes the outward unit normal field at $x \in \Gamma_{1,D}$. Thus, we denote by$$\sigma_1^{D}[-\Delta+c;\mathfrak{B}_1]$$the principal eigenvalue of the problem$$\begin{cases}  -\Delta u + c(x) u=\lambda u& \text{in } D, \\  \mathfrak{B}_1u = 0 & \text{on } \partial D, \end{cases}$$
whose existence and uniqueness are well known (see, for instance, \cite{BNV,Djairo,bookSMP}).
When only the homogeneous Dirichlet condition is considered (i.e., $\Gamma_{1,D} = \emptyset$), we will simply write$$\sigma_1^{D}[-\Delta+c].$$The subsequent result gathers a few properties of this eigenvalue which will be useful throughout this paper.

\begin{lemma}\label{propriedades-local-autovalor}
The following statements hold.
    \begin{enumerate}
        \item[\rm(i)] The principal eigenvalue $\sigma_1^D[-\Delta+c;\mathfrak{B}_1]$ is continuous and strictly increasing with respect to $c \in L^\infty(D)$.
        
        \item[\rm(ii)] If $D_1 \subset D$ is a regular domain such that $\overline{D}_1 \subset D$, then
    $$
    \sigma_1^D[-\Delta+c;\mathfrak{B}_1] <  \sigma_1^{D_1}[-\Delta+c].
    $$ 
    
    \item[\rm(iii)] Defining
    $$
    \mathfrak{P}_{1}:=\{\phi \in W^{2,p}(D): \phi(x)>0 \ \mbox{for all } x \in  D \ \mbox{and } \mathfrak{B}_1 \phi \ge 0 \mbox{ on }\partial D\},
    $$
    then
    \begin{align}\label{supinf-laplaciano}
    \sigma_{1}^{D}[-\Delta+c;\mathfrak{B}_1]= \sup_{\mathfrak{P}_1}\inf_{x \in D}\frac{-\Delta\phi(x)+c(x)\phi(x)}{\phi(x)}.
\end{align}
    \end{enumerate}
\end{lemma}

\begin{proof}
See, for instance, \cite{BNV,Djairo,bookSMP}. 

\end{proof}

\subsection{Local/nonlocal eigenvalue problem} \label{sec2.3}

Let $V:=(v_{1},v_{2}) \in L^{\infty}({\Omega}_\ell) \times C(\overline{\Omega}_{n\ell})$ and consider the following eigenvalue problem:
\begin{equation}
    \label{Pb:autovalor-geral}
    \begin{cases} 
\mathfrak{L}u+V(x)u =\lambda u & \text{in } \Omega,\\ 
\mathfrak{B} u=  0& \text{on } \partial \Omega_{\ell}.
\end{cases}
\end{equation}

The purpose of this section is to recall some fundamental properties of the principal eigenvalue associated with this equation. The results presented here were established in \cite{Peralta2026}. We begin with some definitions.

\begin{definition}
We say that $u \in X$ is a supersolution of $(\mathfrak{{L}}+V, \mathfrak{B}, \Omega)$ if
$$\begin{cases} 
\mathfrak{L}u+V(x)u \ge 0 & \text{in } \Omega, \\ 
\mathfrak{B} u\ge 0 & \text{on } \partial \Omega_\ell.
\end{cases}$$
When one of the above inequalities is strict, we say that $u$ is a strict supersolution.
\end{definition}

 Similarly, we define a subsolution of $(\mathfrak{L}+V, \mathfrak{B}, \Omega)$.

\begin{definition}
We say that $(\mathfrak{L}+V, \mathfrak{B}, \Omega)$ satisfies the Strong Maximum Principle if every supersolution $u\in  X \setminus \{0\}$ of $(\mathfrak{L}+V, \mathfrak{B}, \Omega)$ satisfies $u \gg0$.
\end{definition}

\begin{definition}
    We say that a value $\lambda \in \mathbb{R}$ is an eigenvalue of the operator $(\mathfrak{L}+V,\mathfrak{B},\Omega)$ if the problem 
     \eqref{Pb:autovalor-geral} admits a nontrivial solution in $X_\mathfrak{B}$. Additionally, we say that $\lambda \in \mathbb{R}$ is a principal eigenvalue if it admits a nontrivial solution in $X_\mathfrak{B}$ that does not change sign.
\end{definition}

Next we define
$$
\mathfrak{P}= \{ u \in X:  u(x)>0 \mbox{ for all }x \in \Omega \mbox{ and } \mathfrak{B}u \geq 0 \mbox{ on }\partial\Oe\}.
$$
and 
\begin{equation}
    \lambda_{p}^{\Omega}[\mathfrak{L}+V] := \sup_{\mathfrak{P}} \inf_{x \in \Omega} \frac{(\mathfrak{L}+V)\phi(x)}{\phi(x)}. \label{Def lambdap}
\end{equation}
Whenever there is no risk of confusion, for simplicity, we write $\lambda_{p} := \lambda_{p}^\Omega[\mathfrak{L}+V]$.
Moreover, we denote
$$
m:=\min_{\overline{\Omega}_{n\ell}}v_{2}(x).
$$
Hence,
\begin{equation}
    \label{lambdap<m}\lambda_{p}^{\Omega}[\mathfrak{L}+V] \le m.
\end{equation}

Let us briefly recall some notation. We defined the operators $L_{\ell}, L_{n\ell} : X \to C(\overline{\Omega}_{n\ell})$ as
$$  L_{n\ell} u =-2\int_{\Omega_{n\ell}} J(x-y)u(y)dy \quad \mbox{and} \quad L_{\ell} u= - \int_{\Oe} J(x-y)u(y)dy.  $$
This allows us to rewrite
$$ \mathfrak{L}_2 u = L_{n\ell} u +  L_{\ell} u  \quad \mbox{in } \Onl, $$for every $u \in X$.

The subsequent results can also be derived from \cite{Peralta2026} by setting the functions $B_1$ and $B_2$ in that work as
$$
B_{1} = v_{1} -\int_{\Onl} J(x-y)dy \quad \mbox{and} \quad
B_{2} = v_{2}-2 \int_{\Onl} J(x-y)\,dy - \int_{\Oe} J(x-y)\,dy.
$$

\begin{lemma} \label{Cor-nucleo-trivial}
Suppose that $(\mathfrak{L}+V, \mathfrak{B}, \Omega)$ admits a positive strict supersolution. Then the homogeneous linear problem
\begin{equation}\label{P0}
    \begin{cases}
    \mathfrak{L} u + V(x)u = 0 &\text{in } \Omega, \\
\mathfrak{ B}u = 0 &\text{on } \partial \Oe, \\
\end{cases}
\end{equation}
possesses only the trivial solution $u\equiv0$. That is, $\operatorname{ker}[\mathfrak{L}+V] = \{0\}$, where $\mathfrak{L}+V: X_\mathfrak{B} \to Y$.
\end{lemma}
\begin{proof}
   See \cite[Corollary 2.5]{Peralta2026}.
  
\end{proof}

\begin{proposition}\label{Prop:lambda-av-principal}
If $\lambda \in \mathbb{R}$ is a principal eigenvalue of $(\mathfrak{L}+V,\mathfrak{B},\Omega)$, then $\lambda_p \leq \lambda < m$.
\end{proposition}
\begin{proof}
    See \cite[Proposition 3.4] {Peralta2026}
\end{proof}
\begin{theorem}\label{Th: caso l_p<m}
Assume that $\lambda_p < m$. Then the following statements hold.
\begin{enumerate}
\item[(i)] $\lambda_p$ is a principal eigenvalue of $(\mathfrak{L}+V,\mathfrak{B},\Omega)$ that admits an eigenfunction $\varphi_p \in X_\mathfrak{B}$ such that $\varphi_p \gg0$.

\item[(ii)] $\lambda_p$ is the unique principal eigenvalue of $(\mathfrak{L}+V,\mathfrak{B},\Omega)$.

\item[(iii)] $\lambda_{p}$ is isolated in the spectrum $\sigma(\mathfrak{L}+V)$ and is simple, in the sense that $\operatorname{ker}[\mathfrak{L}+V-\lambda_pI]=\operatorname{span}\{\varphi_p\}$ and the equation
\begin{equation*}
    \begin{cases}
    \mathfrak{L} u + V(x)u - \lambda_p u = \varphi_p  &\text{in } \Omega, \\
\mathfrak{ B}u = 0 &\text{on } \partial \Oe, \\
\end{cases}
\end{equation*}
does not admit a solution $u \in X_\mathfrak{B}$.
\end{enumerate}
\end{theorem}

\begin{proof}
See \cite[Theorem 3.5]{Peralta2026}.
\end{proof}

The following result provides sufficient conditions for $\lambda_{p}^{\Omega}[\mathfrak{L}+V]<m$.

\begin{theorem}\label{Condi sufi}
If there exists a measurable set $A \subset \Onl$ such that
\begin{align}\label{lambdap<m-cond suficiente}
   \frac{1}{2} < \inf_{y \in A} \int_A \frac{J(x-y)}{v_{2}(x)-m} dx \le \infty, 
\end{align}
then $\lambda_{p}<m$.
\end{theorem}

\begin{proof}
See \cite[Theorem 3.7]{Peralta2026}.
\end{proof}

Moreover, $\lambda_{p}$ satisfies the following properties.

\begin{proposition}\label{Propriedades}
The following properties hold:
\begin{enumerate}
    \item[(i)] If $c \in \mathbb{R}$, then
    $$
    \lambda_{p}^{\Omega}[\mathfrak{L}+V+c] = \lambda_{p}^{\Omega}[\mathfrak{L}+V] + c.
    $$

    \item[(ii)] Let $v_1, w_1 \in L^{\infty}(\Omega_{\ell})$ and $v_2, w_2 \in C(\overline{\Omega}_{n\ell})$, and define $V := (v_{1}, v_{2})$ and $W := (w_{1}, w_{2})$. If $V \le W$ in $\Omega$, then
    $$
    \lambda_{p}^{\Omega}[\mathfrak{L}+V] \le \lambda_{p}^{\Omega}[\mathfrak{L}+W].
    $$

    \item[(iii)] Let $v_{1,k}, v_1 \in L^{\infty}(\Omega_{\ell})$ and $v_{2,k}, v_2 \in C(\overline{\Omega}_{n\ell})$, $k \in \mathbb{N}$, and define $V_k:= (v_{1,k}, v_{2,k})$ and $V := (v_1, v_2)$. If $V_{k} \to V$ in $L^{\infty}(\Omega_{\ell
    }) \times C(\overline{\Omega}_{n\ell})$, then
    $$
    \lambda_{p}^{\Omega}[\mathfrak{L}+V_k] \to \lambda_{p}^{\Omega}[\mathfrak{L}+V].
    $$
\end{enumerate}
\end{proposition}

\begin{proof}
See \cite[Proposition 4.2]{Peralta2026}.
\end{proof}

\begin{proposition}\label{autovalor com supersoluçao}
    Let $$
\mathcal{A} := \Big\{ \lambda \in \mathbb{R} : \exists \, \phi \in X, \, \phi(x) > 0, \, x \in \Omega,
     \text{ supersolution of } (\mathfrak{L} + V - \lambda, \mathfrak{B}, \Omega) \Big\}.
$$
Then, $$\mu_{p}^{\Omega}[\mathfrak{L}+V]:=\sup \mathcal{A}= \lambda_{p}^{\Omega}[\mathfrak{L}+V].$$
\end{proposition}
\begin{proof}
    See \cite[Proposition 4.1]{Peralta2026}.
\end{proof}

\section{Strict monotonicity properties}\label{sec3}

The main goal of this section is to provide sufficient conditions ensuring that some monotonicity properties of $\lambda_p^\Omega[\mathfrak{L}+V]$ are strict. Our results show that, if one of the principal eigenvalues involved admits an eigenfunction, then the corresponding monotonicity is strict.

\begin{figure}
    \centering
 
\tikzset{
pattern size/.store in=\mcSize, 
pattern size = 5pt,
pattern thickness/.store in=\mcThickness, 
pattern thickness = 0.3pt,
pattern radius/.store in=\mcRadius, 
pattern radius = 1pt}
\makeatletter
\pgfutil@ifundefined{pgf@pattern@name@_ejjf6l9as lines}{
\pgfdeclarepatternformonly[\mcThickness,\mcSize]{_ejjf6l9as}
{\pgfqpoint{0pt}{0pt}}
{\pgfpoint{\mcSize+\mcThickness}{\mcSize+\mcThickness}}
{\pgfpoint{\mcSize}{\mcSize}}
{\pgfsetcolor{\tikz@pattern@color}
\pgfsetlinewidth{\mcThickness}
\pgfpathmoveto{\pgfpointorigin}
\pgfpathlineto{\pgfpoint{\mcSize}{0}}
\pgfusepath{stroke}}}
\makeatother

 
\tikzset{
pattern size/.store in=\mcSize, 
pattern size = 5pt,
pattern thickness/.store in=\mcThickness, 
pattern thickness = 0.3pt,
pattern radius/.store in=\mcRadius, 
pattern radius = 1pt}
\makeatletter
\pgfutil@ifundefined{pgf@pattern@name@_ce0ds7zyy}{
\pgfdeclarepatternformonly[\mcThickness,\mcSize]{_ce0ds7zyy}
{\pgfqpoint{-\mcThickness}{-\mcThickness}}
{\pgfpoint{\mcSize}{\mcSize}}
{\pgfpoint{\mcSize}{\mcSize}}
{
\pgfsetcolor{\tikz@pattern@color}
\pgfsetlinewidth{\mcThickness}
\pgfpathmoveto{\pgfpointorigin}
\pgfpathlineto{\pgfpoint{0}{\mcSize}}
\pgfusepath{stroke}
}}
\makeatother
\tikzset{every picture/.style={line width=0.75pt}} 

\begin{tikzpicture}[x=0.75pt,y=0.75pt,yscale=-1,xscale=1]

\draw  [fill={rgb, 255:red, 128; green, 128; blue, 128 }  ,fill opacity=1 ] (190.81,138.14) .. controls (190.81,83.94) and (241.99,40) .. (305.11,40) .. controls (368.24,40) and (419.41,83.94) .. (419.41,138.14) .. controls (419.41,192.34) and (368.24,236.28) .. (305.11,236.28) .. controls (241.99,236.28) and (190.81,192.34) .. (190.81,138.14) -- cycle ;
\draw  [color={rgb, 255:red, 0; green, 0; blue, 0 }  ,draw opacity=1 ][pattern=_ejjf6l9as,pattern size=6pt,pattern thickness=0.75pt,pattern radius=0pt, pattern color={rgb, 255:red, 0; green, 0; blue, 0}][line width=0.75]  (223.61,137.17) .. controls (223.61,97.92) and (261.13,66.11) .. (307.42,66.11) .. controls (353.72,66.11) and (391.24,97.92) .. (391.24,137.17) .. controls (391.24,176.42) and (353.72,208.24) .. (307.42,208.24) .. controls (261.13,208.24) and (223.61,176.42) .. (223.61,137.17) -- cycle ;
\draw    (341.48,170.36) -- (367.24,203.53) -- (395.91,240.45) ;
\draw  [pattern=_ce0ds7zyy,pattern size=6pt,pattern thickness=0.75pt,pattern radius=0pt, pattern color={rgb, 255:red, 0; green, 0; blue, 0}][line width=1.5]  (250.68,138.14) .. controls (250.68,114.91) and (275.05,96.08) .. (305.11,96.08) .. controls (335.17,96.08) and (359.54,114.91) .. (359.54,138.14) .. controls (359.54,161.37) and (335.17,180.2) .. (305.11,180.2) .. controls (275.05,180.2) and (250.68,161.37) .. (250.68,138.14) -- cycle ;
\draw  [color={rgb, 255:red, 0; green, 0; blue, 0 }  ,draw opacity=1 ][fill={rgb, 255:red, 128; green, 128; blue, 128 }  ,fill opacity=1 ][line width=0.75]  (277.9,137.17) .. controls (277.9,125.56) and (291.12,116.14) .. (307.42,116.14) .. controls (323.73,116.14) and (336.95,125.56) .. (336.95,137.17) .. controls (336.95,148.79) and (323.73,158.2) .. (307.42,158.2) .. controls (291.12,158.2) and (277.9,148.79) .. (277.9,137.17) -- cycle ;
\draw    (419,130) -- (429,130) ;

\draw (378.52,46.27) node [anchor=north west][inner sep=0.75pt]    {$\Omega $};
\draw (272.8,153.18) node [anchor=north west][inner sep=0.75pt]  [font=\large,color={rgb, 255:red, 0; green, 0; blue, 0 }  ,opacity=1 ]  {$\Omega _{n\ell *}$};
\draw (299.27,127.09) node [anchor=north west][inner sep=0.75pt]    {$\Omega _{n\ell }$};
\draw (281,69.92) node [anchor=north west][inner sep=0.75pt]  [color={rgb, 255:red, 0; green, 0; blue, 0 }  ,opacity=1 ]  {$\Omega _{\ell *}$};
\draw (197.92,104.41) node [anchor=north west][inner sep=0.75pt]    {$\Omega {_{\ell }}$};
\draw (397,231.4) node [anchor=north west][inner sep=0.75pt]    {$\Gamma $};
\draw (429,119.4) node [anchor=north west][inner sep=0.75pt]    {$\Gamma _{0}$};

\end{tikzpicture}
    \caption{Possible configuration for the subdomains.}
    \label{fig: sub Dominio}
\end{figure}
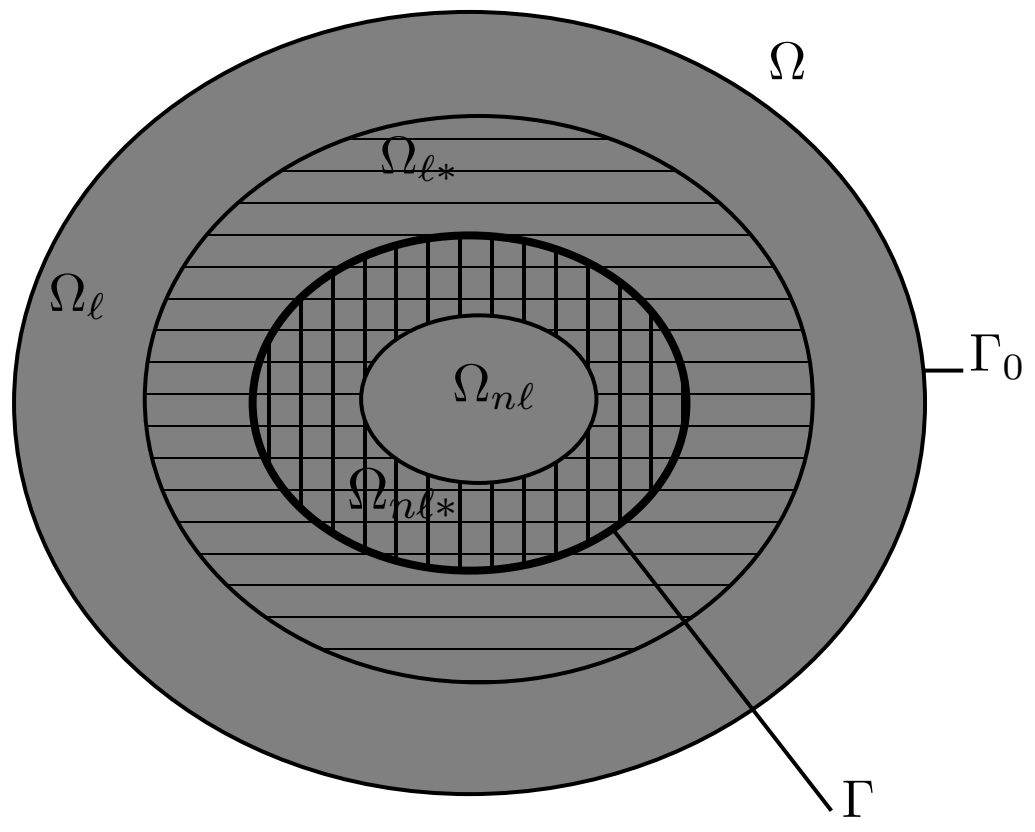

\begin{proposition}\label{monotonia-estrita-peso}
        Let $V=(v_{1}(x),v_{2}(x)),W:=(w_{1}(x),w_{2}(x)) \in L^\infty({\Omega}_{\ell}) \times C(\overline{\Omega}_{n\ell})$ be such that $V < W$. If $\lambda_{p}^{\Omega}[\mathfrak{L}+V]$ or $\lambda_{p}^{\Omega}[\mathfrak{L}+W]$ admits a positive eigenfunction, then $$\lambda_{p}^{\Omega}[\mathfrak{L}+V] < \lambda_{p}^{\Omega}[\mathfrak{L}+W].$$
\end{proposition} 
\begin{proof}
    Since $V < W$, by Proposition \ref{Propriedades}, it follows that $\lambda_{p}^{\Omega}[\mathfrak{L}+V] \le \lambda_{p}^{\Omega}[\mathfrak{L}+W].$
   Suppose for the sake of contradiction that 
   \begin{equation}\label{lambdaV=lambdaW}
       \lambda_{p}^{\Omega}[\mathfrak{L}+V]=\lambda_{p}^{\Omega}[\mathfrak{L}+W].
   \end{equation}
For each $\epsilon \in [0,1]$, define $$\lambda(\epsilon):=\lambda_{p}^{\Omega}[\mathfrak{L}+ \epsilon W  +(1-\epsilon)V] \quad \mbox{and }\quad m(\epsilon):=\min_{\overline{\Omega}_{n\ell}}(\epsilon w_{2}(x)+(1-\epsilon)v_{2}(x)).$$
Since 
$$V \le \epsilon W  +(1-\epsilon)V \le W \quad \mbox{in }\Omega, \quad \forall  \epsilon \in [0,1],$$
using Proposition \ref{Propriedades} again, we have that
   $$\lambda_{p}^{\Omega}[\mathfrak{L}+V] \le \lambda(\epsilon) \le \lambda_{p}^{\Omega}[\mathfrak{L}+W] \quad \forall \epsilon \in [0,1].$$
Therefore, \eqref{lambdaV=lambdaW} implies that 
\begin{equation}
    \lambda(\epsilon)=\lambda_{p}^{\Omega}[\mathfrak{L}+W] \quad \forall \epsilon \in [0,1].
\end{equation}
Now suppose that $\lambda_p^\Omega[\mathfrak{L}+W]$ admits an eigenfunction $\varphi_W \in X_\mathfrak{B}$ such that $\varphi_W \gg0$. By Proposition \ref{Prop:lambda-av-principal}, we have 
$$\lambda(1)=  \lambda_{p}^\Omega[\mathfrak{L}+W]< \min_{\overline{\Omega}_{n\ell}}(w_{2}(x)) = m(1).$$
By continuity, there exists $\epsilon_{0}>0$ such that   $$\lambda_p^\Omega[\mathfrak{L}+W]=\lambda(\epsilon)<m(\epsilon), \quad \forall \epsilon \in (\epsilon_{0},1].$$
Using Theorem \ref{Th: caso l_p<m}, for each $\epsilon \in (\epsilon_0,1)$, there exists an eigenfunction $\varphi_{\epsilon} \in X_{\mathfrak{B}}$ such that $\varphi_{\epsilon} \gg0$ satisfying 
 \begin{equation*}
    \begin{cases}
   (\mathfrak{L}+(\epsilon W+(1-\epsilon)V))\varphi_{\epsilon}=\lambda(\epsilon) \varphi_{\epsilon}=\lambda_{p}^{\Omega}[\mathfrak{L}+W]\varphi_{\epsilon}  &\text{in } \Omega, \\
\mathfrak{ B}\varphi_{\epsilon}= 0 &\text{on } \partial \Oe. \\
\end{cases}
\end{equation*}
Since $\epsilon <1$, we have that $$(\mathfrak{L}+W-\lambda_{p}^{\Omega}[\mathfrak{L}+W])\varphi_{\epsilon}>(\mathfrak{L}+(\epsilon W+(1-\epsilon)V))\varphi_{\epsilon}-\lambda_{p}^{\Omega}[\mathfrak{L}+W]\varphi_{\epsilon}=0 \quad\mbox{in }\Omega. $$
 Therefore, we deduce that $\varphi_{\epsilon} \gg0$ is a strict supersolution for the operator $(\mathfrak{L}+W-\lambda_{p}^{\Omega}[\mathfrak{L}+W],\mathfrak{B}, \Omega)$. By Lemma \ref{Cor-nucleo-trivial}, we obtain that the problem 
 \begin{equation*}
    \begin{cases}
    \mathfrak{L} u + W(x)u - \lambda_p^{\Omega}[\mathfrak{L}+W]u = 0 &\text{in } \Omega, \\
\mathfrak{ B}u = 0 &\text{on } \partial \Oe, \\
\end{cases}
\end{equation*}
admits only the trivial solution $u=0$, which contradicts the assumption that $\varphi_{W} \gg0$ is an eigenfunction associated with $\lambda_p^\Omega[\mathfrak{L}+W]$. This concludes the proof for this case.

Now suppose that $\lambda_p^\Omega[\mathfrak{L}+V]$ admits an eigenfunction $\varphi_V \in X_\mathfrak{B}$ such that $\varphi_V\gg0$. Arguing as in the previous case, we have that $\lambda(0)=\lambda_p^\Omega[\mathfrak{L}+V]<\min_{\overline{\Omega}_{n\ell}}v_2(x)$ and there exists $\epsilon_0>0$ such that
$$
\lambda(\epsilon
)= \lambda_p^\Omega[\mathfrak{L}+V] <m(\epsilon) \quad \forall \epsilon\in [0,\epsilon_0).$$
Thus, for each $\epsilon \in (0,\epsilon_0)$, there exists an eigenfunction $\varphi_\epsilon\in X_\mathfrak{B}$ with $\varphi_\epsilon\gg0$ satisfying
 \begin{equation*}
    \begin{cases}
   (\mathfrak{L}+(\epsilon W+(1-\epsilon)V)-\lambda_{p}^{\Omega}[\mathfrak{L}+V])\varphi_{\epsilon} =0 &\text{in } \Omega, \\
\mathfrak{ B}\varphi_{\epsilon}= 0 &\text{on } \partial \Oe. \\
\end{cases}
\end{equation*}
However, since $\epsilon>0$ and $V<W$, a direct calculation shows that $\varphi_V\gg0$ is a strict positive supersolution of $(\mathfrak{L}+(\epsilon W+(1-\epsilon)V)-\lambda_{p}^{\Omega}[\mathfrak{L}+V],\mathfrak{B},\Omega)$. Thus, Lemma~\ref{Cor-nucleo-trivial} implies that the above equation admits only the trivial solution, which is a contradiction. This completes the proof.
\end{proof}
\begin{remark}
    If $\lambda_{p}^{\Omega}[\mathfrak{L}+W]$ and $\lambda_{p}^{\Omega}[\mathfrak{L}+V]$ do not admit eigenfunctions, we combine  Theorem \ref{Th: caso l_p<m} with \eqref{lambdap<m} to obtain$$\lambda_{p}^{\Omega}[\mathfrak{L}+W]=\min_{\overline{\Omega}_{n\ell}}w_{2}(x) \quad \text{and} \quad \lambda_{p}^{\Omega}[\mathfrak{L}+V]=\min_{\overline{\Omega}_{n\ell}}v_{2}(x).$$
    Consequently, strict monotonicity holds if and only if     $$\min_{\overline{\Omega}_{n\ell}}v_{2}(x)<\min_{\overline{\Omega}_{n\ell}}w_{2}(x).$$
\end{remark}
Subsequently, we investigate when the monotonicity of $\lambda_p^\Omega[\mathfrak{L}+W]$ with respect to the domain is strict. To this end, let $\Omega_{n\ell*} \subset \Onl$ and $\Omega_{\ell*} \subset \Oe$ be regular subdomains  with $\operatorname{dist}(\Omega_{\ell*},\Omega_{n\ell*})<\delta$ such that either $\Gamma \subset \partial \Omega_{\ell*} \cap \partial \Omega_{n\ell*}   $, or  $\Gamma \cap \partial \Omega_{\ell*} \cap \partial \Omega_{n\ell*} = \emptyset$  (see Figure \ref{fig: sub Dominio}). Setting $\Omega_*= \operatorname{int}(\overline{\Omega_{\ell*} \cup \Omega_{n\ell*}}) \subset \Omega$, we define the sets
\[
    X_*= \Big\{u \in L^2(\Omega_*) : u|_{\Omega_{\ell*}} \in W^{2,p}(\Omega_{\ell*}),\, p > N, \text{ and }  u|_{\overline{\Omega}_{n\ell*}} \in 
    C(\overline{\Omega}_{n\ell*}) \Big\},
\]
$$\mathfrak{P}_*= \Big\{ u \in X_*:  u(x)>0 \mbox{ for all }x \in \Omega_{*} \mbox{ and } \mathfrak{B}_*u \geq 0 \mbox{ on }\partial\Omega_{\ell*}\Big\}$$
and the operators
\[
 \mathfrak{L}_*u:= \begin{cases} -\Delta u - \displaystyle\int_{\Omega_{n\ell*}}  J(x-y)u(y)dy    & \text{in } \Omega_{\ell*}, \\ -\displaystyle2\int_{\Omega_{n\ell*}} J(x-y)u(y)dy - \int_{\Omega_{\ell*}} J(x-y)u(y)dy  & \text{in } \Omega_{n \ell*}, \end{cases}
\]
and
\begin{equation*}
\mathfrak{B}_* u := \begin{cases} 
\displaystyle \frac{\partial u}{\partial \eta} & \text{on } \Gamma,\\
u & \text{on } \partial \Omega_{\ell*} \setminus \Gamma.
\end{cases}
\end{equation*}

\begin{proposition}\label{monotonia-estrita-dominio}
    Let $\Omega_{*} \subset \Omega$ be a proper subdomain as above. Then 
    $$\lambda_{p}^{\Omega}[\mathfrak{L}+V] \leq \lambda_{p}^{\Omega_{*}}[\mathfrak{L}_{*}+V].$$
    Moreover, if $\lambda_{p}^{\Omega}[\mathfrak{L}+V]$ or $\lambda_{p}^{\Omega_{*}}[\mathfrak{L}_{*}+V]$ admits a positive eigenfunction, then the inequality is strict.   
\end{proposition}
\begin{proof}
    Under the above assumptions, it  follows from \cite[Corollary 4.3]{Peralta2026} that 
    $$\lambda_{p}^\Omega := \lambda_{p}^{\Omega}[\mathfrak{L}+V] \le \lambda_{p}^{\Omega_{*}}[\mathfrak{L}_{*}+V] := \lambda_{p}^{\Omega_{*}}.$$
    In what follows, we prove that the inequality is strict. To this end, we argue by contradiction, assuming that 
$$\lambda_{p}^{\Omega} = \lambda_{p}^{\Omega_*}.$$
    Let us consider the case where $\lambda_{p}^{\Omega}$ admits a positive eigenfunction $\varphi \in X_{\mathfrak{B}}$, with $\varphi \gg 0$. Combining $\Omega_{n\ell*} \subset \Onl$ with Proposition \ref{Prop:lambda-av-principal}, we obtain that 
    $$\lambda_{p}^{\Omega_{*}} = \lambda_{p}^{\Omega} < \min_{\overline{\Omega}_{n\ell}}v_{2}(x) \le \min_{\overline{\Omega}_{n\ell*}}v_{2}(x).$$
    Consequently, by Theorem \ref{Th: caso l_p<m}, there exists a positive eigenfunction $\varphi_{*} \in X_{\mathfrak{B}_{*}}$, $\varphi_{*} \gg 0$, associated with $\lambda_{p}^{\Omega_{*}}$ satisfying
    \begin{equation*}
    \begin{cases}
    \mathfrak{L}_{*} \varphi_{*} + V(x)\varphi_{*} - \lambda_p^{\Omega_*}\varphi_{*} = 0 &\text{in } \Omega_{*}, \\
    \mathfrak{B}_*\varphi_{*} = 0 &\text{on } \partial \Omega_{\ell*}.
    \end{cases}
    \end{equation*}
    On the other hand, since $\Omega_{\ell*} \subset \Oe$ and $\Omega_{n\ell*} \subset \Onl$, we find that 
    \begin{align*}
        - \int_{\Omega_{n\ell*}} J(x-y)\varphi(y)dy &\ge -\int_{\Omega_{n\ell}} J(x-y)\varphi(y)dy \quad \text{for } x \in \Omega_{n\ell*} \text{ and } x \in \Omega_{\ell*}, \\ 
        -\int_{\Omega_{\ell*}} J(x-y)\varphi(y)dy &\ge -\int_{\Omega_{\ell}} J(x-y)\varphi(y)dy \quad \text{for } x \in \Omega_{n\ell*}.
    \end{align*}
    Moreover, since the inclusion $\Omega_* \subset \Omega$ is proper, one of the inequalities above will be strict for $x \in \Omega_*$ such that $\operatorname{dist}(x, \Omega \setminus\Omega_*) < \delta$. Consequently, a direct computation shows that $\varphi|_{\Omega_*}$ satisfies
    \begin{equation*}
    \begin{cases}
    \mathfrak{L}_{*} \varphi + V(x)\varphi - \lambda_p^{\Omega_*}[\mathfrak{L}_{*}+V]\varphi > 0 &\text{in } \Omega_{*}, \\
    \mathfrak{B}_*\varphi > 0 &\text{on } \partial \Omega_{\ell*}.
    \end{cases}
    \end{equation*}
    That is, $\varphi|_{\Omega_{*}}$ is a strict positive supersolution for the operator $(\mathfrak{L}_{*}+V-\lambda_{p}^{\Omega_{*}}[\mathfrak{L}_*+V], \mathfrak{B}_*, \Omega_{*})$. By Lemma~\ref{Cor-nucleo-trivial}, we infer that the problem  
    \begin{equation*}
    \begin{cases}
    \mathfrak{L}_{*} u + V(x)u - \lambda_p^{\Omega_*}[\mathfrak{L}_{*}+V]u = 0 &\text{in } \Omega_{*}, \\
    \mathfrak{B}_{*}u = 0 &\text{on } \partial \Omega_{\ell*},
    \end{cases}
    \end{equation*}
    admits only the trivial solution $u = 0$, which is a contradiction, given that $\varphi_{*} \gg 0$ is a solution to this equation. 

Now we deal with the case where $\lambda_{p}^{\Omega_{*}}$ admits a positive eigenfunction $\varphi_{*} \in X_{\mathfrak{B}_{*}}$, $\varphi_* \gg 0$. Thus, by Proposition \ref{Prop:lambda-av-principal}, we obtain 
    \begin{equation}
        \label{lp-lp*}\lambda_{p}^{\Omega} = \lambda_{p}^{\Omega_{*}} < \min_{\overline{\Omega}_{n\ell*}}v_{2}(x).
    \end{equation}
    Note that $\varphi_{*}$ solves the following problem 
    \begin{equation}
    \begin{cases}
    \mathfrak{L}_{*} u + V(x)u - \lambda_p^{\Omega}[\mathfrak{L}+V]u = 0 &\text{in } \Omega_{*}, \\
    \mathfrak{B}_{*}u = 0 &\text{on } \partial \Omega_{\ell*}.
    \end{cases} \label{eq:operador-estrela}
    \end{equation}
    Let us consider the following possibilities:

    \textbf{Case $\Omega_{n\ell*} = \Onl$}. It follows from \eqref{lp-lp*} that 
    $$\lambda_{p}^{\Omega} < \min_{\overline{\Omega}_{n\ell}}v_{2}(x).$$
    By Theorem \ref{Th: caso l_p<m}, there exists an eigenfunction $\varphi \in X_{\mathfrak{B}}$, $\varphi \gg 0$ associated with $\lambda_{p}^{\Omega}$. 
    Since $\Omega_{\ell*} \subset \Oe$, we can argue as before to conclude that $\varphi|_{\Omega_*}$ is a strict positive supersolution for the operator $(\mathfrak{L}_{*}+V-\lambda_{p}^\Omega, \mathfrak{B}_*, \Omega_*)$. By Lemma~\ref{Cor-nucleo-trivial}, \eqref{eq:operador-estrela} admits only the trivial solution $u \equiv 0$, contradicting the fact that $\varphi_{*} \gg 0$ is a solution to this equation.

    \textbf{Case $\Omega_{n\ell*} \varsubsetneq \Onl$}.
    Let $\epsilon > 0$ be a constant, and define 
    $$ \Omega_{n\ell}^{\epsilon} := \{x \in \Onl : \operatorname{dist}(x,\Omega_{n\ell*}) < \epsilon\} \quad \text{and} \quad \Omega^\epsilon := \Omega_{\ell*} \cup \Omega_{n\ell}^\epsilon.$$
    Consequently, $$\Omega_{*} \varsubsetneq \Omega^\epsilon \subset \Omega.$$
    For sufficiently small $\epsilon>0$, $\Omega^\epsilon$ is an admissible configuration for the coupled local/nonlocal models. Thus,  the monotonicity already established implies that
    \begin{align}
        \lambda_{p}^{\Omega} \le \lambda_{p}^{\epsilon} \le \lambda_{p}^{\Omega_*} = \lambda_{p}^{\Omega} < \min_{\overline{\Omega}_{n\ell*}}v_{2}(x), \label{eq:estimativas-dominios}
    \end{align}
    where $\lambda_p^\epsilon := \lambda_p^{\Omega^\epsilon}[\mathfrak{L}_\epsilon+V]$ and $\mathfrak{L}_\epsilon$ is defined analogously to $\mathfrak{L}_*$, with the integrals over $\Omega_{n\ell*}$ replaced by $\Omega_{n\ell}^{\epsilon}$. 
    Since 
    $$\min_{\overline{\Omega}_{n\ell}^\epsilon}v_{2}(x) \to \min_{\overline{\Omega}_{n\ell*}}v_{2}(x) \quad \text{as } \epsilon \to 0^+,$$
    we deduce from \eqref{eq:estimativas-dominios} that for $\epsilon$ sufficiently small, $$\lambda_{p}^{{\epsilon}} < \min_{\overline{\Omega}_{n\ell}^\epsilon}v_{2}(x).$$
    Applying Theorem \ref{Th: caso l_p<m} again, there exists $\varphi_{\epsilon} \in X_{\mathfrak{B}}$, $\varphi_\epsilon \gg 0$, an eigenfunction of $\lambda_{p}^{\Omega^{\epsilon}} = \lambda_p^{\Omega}$ satisfying
    \begin{equation}\label{varphi_epsilon}
    \begin{cases}
    \mathfrak{L}_\epsilon \varphi_{\epsilon} + V(x)\varphi_{\epsilon} - \lambda_p^{\Omega}\varphi_{\epsilon} = 0 &\text{in } \Omega^{\epsilon}, \\
    \mathfrak{B}_*\varphi_{\epsilon} = 0 &\text{on } \partial \Omega_{\ell*}.
    \end{cases}
    \end{equation}
    On the other hand, since $\Omega_{n\ell*} \varsubsetneq \Omega_{n\ell}^{\epsilon}$, it follows that 
    $$-2\int_{\Omega_{n\ell*}} J(x-y)\varphi_{\epsilon}(y)dy \geq -2\int_{\Omega_{n\ell}^{\epsilon}} J(x-y)\varphi_{\epsilon}(y)dy \quad \text{in } \Omega_{n\ell*},$$
    with the inequality being strict for all $x \in \Omega_{n\ell*}$ such that $\operatorname{dist}(x, \Omega_{n\ell}^{\epsilon} \setminus \Omega_{n\ell{*}}) < \delta$. An analogous inequality holds for the other integrals. Consequently,
    $$
    \mathfrak{L}_* \varphi_\epsilon > \mathfrak{L}_\epsilon \varphi_\epsilon \quad \text{in }\Omega_*.
    $$
    Combining this inequality with \eqref{varphi_epsilon}, we conclude that $\varphi_\epsilon|_{\Omega_*}$ is a strict positive supersolution for $(\mathfrak{L}_*+V-\lambda_{p}^{\Omega}, \mathfrak{B}_*, \Omega_*).$ By Lemma~\ref{Cor-nucleo-trivial}, problem \eqref{eq:operador-estrela} admits only the solution $u \equiv 0$, which is a contradiction, given that $\varphi_{*} \gg 0$ is a solution to this equation.
    This concludes the proof.
\end{proof}

We conclude this section with a result relating the generalized eigenvalues introduced in Section~\ref{Secao2}.

\begin{proposition}\label{estimativas}
The following inequalities hold:
\begin{enumerate} 
    \item [\rm(i)] $\lambda_{p}^\Omega[\mathfrak{L}+V] < \sigma_{1}^{\Omega_{\ell}}[-\Delta+v_{1};\mathfrak{B}]$. 
    \item [\rm(ii)] $\lambda_{p}^\Omega[\mathfrak{L}+V]\le \lambda_{p}^{\Onl}[L_{n\ell}+v_{2}]$. Moreover, if $\lambda_{p}^\Omega[\mathfrak{L}+V]$ or $ \lambda_{p}^{\Onl}[L_{n\ell}+v_{2}]$ admits an eigenfunction, then the inequality is strict.
\end{enumerate}

\end{proposition}

\begin{proof}
We first prove $(i)$. Let $u \in \mathfrak{P}$. For every $x \in \Oe$, we have
$$
\frac{(\mathfrak{L}+V)u(x)}{u(x)}= \frac{-\Delta u(x)-\int_{\Onl}J(x-y)u(y)dy+v_{1}(x)u(x)}{u(x)} \le\frac{-\Delta u(x)+v_1(x)u(x)}{u(x)}.
$$
Since $\Oe \subset \Omega$, \eqref{supinf-laplaciano} implies that
$$
\inf_{\Omega}\frac{(\mathfrak{L}+V)u(x)}{u(x)}\le \inf_{\Oe}\frac{-\Delta u(x)+v_1(x)u(x)}{u(x)} \le \sigma_{1}^{\Oe}[-\Delta+v_{1};\mathfrak{B}].
$$
Hence, from definition \eqref{Def lambdap}, it follows that
$$
\lambda_{p}^{\Omega}[\mathfrak{L}+V] \le \sigma_{1}^{\Oe}[-\Delta+v_{1};\mathfrak{B}].
$$
Now, we will show that this inequality is strict. Suppose, for the sake of contradiction, that
$$
\lambda_{p}^{\Omega}[\mathfrak{L}+V] = \sigma_{1}^{\Oe}[-\Delta+v_{1};\mathfrak{B}]=: \sigma_{1}^{\Oe}.
$$
It follows from \eqref{lambdap<m} that $$\sigma_1^{\Omega_\ell}=\lambda_{p}^\Omega[\mathfrak{L}+V] \le \min_{\overline{\Omega}_{n\ell}} v_2(x).$$
Let us analyze the two possible cases:

\textbf{Case 1:} $\sigma_{1}^{\Oe} < \min_{\overline{\Omega}_{n\ell}} v_2(x)$.

Then  $\lambda_{p}^{\Omega}[\mathfrak{L}+V] < \min_{\overline{\Omega}_{n\ell}} v_2(x)$. By Theorem \ref{Th: caso l_p<m}, there exists an eigenfunction $\varphi \in X_\mathfrak{B}$, $\varphi \gg 0$, associated with $\lambda_{p}^{\Omega}$ such that $|\varphi|_{2}=1$. Since $\operatorname{dist}(\Oe, \Omega_{n\ell}) < \delta$, we deduce that $\varphi|_{\Oe}$ satisfies
\begin{align*}
\begin{cases}
-\Delta \varphi+v_{1}(x) \varphi-\sigma_{1}^{\Oe}\varphi= \int_{\Onl} J(x-y) \varphi(y)dy>0 \quad &\mbox{in } \Oe,
\\
\mathfrak{B}\varphi=0 \quad &\mbox{on } \partial \Oe.
\end{cases}
\end{align*}
Therefore, $\varphi|_{\Oe}$ is a strict positive supersolution for the local operator $(-\Delta+v_{1}(x)-\sigma_{1}^{\Oe},\mathfrak{B}, \Oe)$. By the characterization of the maximum principle (see, for instance, \cite{bookSMP}), we infer that
$$
0 = \sigma_{1}^{\Oe}[-\Delta+v_1; \mathfrak{B}]-\sigma_{1}^{\Oe} = \sigma_{1}^{\Oe}[-\Delta+v_{1}-\sigma_{1}^{\Oe}; \mathfrak{B}]>0,
$$
which is a contradiction.

\textbf{Case 2:} $\sigma_{1}^{\Oe} = \min_{\overline{\Omega}_{n\ell}} v_2(x)$. 

Thus $\lambda_{p}^{\Omega}[\mathfrak{L}+V] = \min_{\overline{\Omega}_{n\ell}} v_2(x)$. Let $\epsilon>0$ and define $V_{\epsilon}:=(v_{1}(x),v_{2}(x)+\epsilon).$ Since $V \le V_{\epsilon}$, we can combine Proposition \ref{Propriedades} (ii) with what we have already proved to conclude that
$$\lambda_{p}^{\Omega}[\mathfrak{L}+V] \le \lambda_{p}^{\Omega}[\mathfrak{L}+V_{\epsilon}]\le \sigma_{1}^{\Oe}=\min_{\overline{\Omega}_{n\ell}} v_{2}(x)<\min_{\overline{\Omega}_{n\ell}} v_{2}(x)+\epsilon=:m_{\epsilon}.$$
Consequently, $\lambda_p^\Omega[\mathfrak{L}+V_\epsilon]$ admits a positive eigenfunction, and Proposition \ref{monotonia-estrita-peso} implies that
$$
\lambda_p^\Omega[\mathfrak{L}+V] < \lambda_p^\Omega[\mathfrak{L}+V_\epsilon] \leq \sigma_1^{\Omega_\ell},
$$
contradicting the initial assumption that $\lambda_{p}^{\Omega}[\mathfrak{L}+V]= \sigma_{1}^{\Oe}.$ This concludes the proof of item (i).

We now prove (ii). First, we show that the inequality holds. Given $u \in \mathfrak{P}$, we also have
$$
\frac{(\mathfrak{L}+V)u(x)}{u(x)}= \frac{L_{n\ell}u(x)+L_{\ell}u(x)+v_{2}(x)u(x)}{u(x)} \le\frac{L_{n\ell} u(x)+v_2(x)u(x)}{u(x)} \quad \mbox{in }\Onl.
$$
Therefore, from  \eqref{autovalor nonlocal} and \eqref{Def lambdap}, it follows that
$$
\lambda_{p}^{\Omega}[\mathfrak{L}+V] \le\lambda_{p}^{\Onl}[L_{n\ell}+v_{2}].
$$
We will prove that the strict inequality holds. To this end, suppose for the sake of contradiction that
$$\lambda_{p}^{\Omega}[\mathfrak{L}+V]=\lambda_{p}^{\Onl}[L_{n\ell}+v_{2}]=: \lambda_{p}^{\Onl}.$$
If $\lambda_{p}^{\Omega}[\mathfrak{L}+V]$ admits a positive eigenfunction $\varphi \in X_{\mathfrak{B}}$, $\varphi \gg 0$, it follows from Proposition \ref{Prop:lambda-av-principal} that
    $$\lambda_{p}^{\Onl}=\lambda_{p}^{\Omega}[\mathfrak{L}+V]<\min_{\overline{\Omega}_{n\ell}}(v_{2}(x)).$$ 
Applying Lemma \ref{propriedades-NL-autovalor} (iii), there exists a positive eigenfunction $\varphi_{n\ell} \in C(\overline{\Omega}_{n\ell})$ satisfying
$$L_{n\ell} \varphi_{n\ell}+v_{2}(x) \varphi_{n\ell}-\lambda_{p}^{\Onl}\varphi_{n\ell}=0 \quad  \mbox{in } \overline{\Omega}_{n\ell}.$$
On the other hand, since $\operatorname{dist}(\Oe, \Onl)<\delta$, it follows from assumption \ref{J1} that
$$L_{\ell}\varphi(x) = -\int_{\Oe} J(x-y) \varphi(y) dy<0$$
for all $x \in \Onl$ such that $\operatorname{dist}(x,\Oe)<\delta$. Thus, $\varphi|_{\Onl}$ satisfies
\begin{align*}
    L_{n\ell} \varphi+v_{2}(x) \varphi-\lambda_{p}^{\Omega_{n\ell}}\varphi > L_{n\ell} \varphi+ L_{\ell}\varphi+v_{2}(x) \varphi-\lambda_{p}^{\Omega_{n\ell}}\varphi=0 \quad \mbox{in } \overline{\Omega}_{n\ell
    }.
\end{align*}
Therefore, we infer that $\varphi|_{\Onl}$ is a strict positive supersolution for the operator $L_{n\ell}+v_{2}(x)-\lambda_{p}^{\Onl}$ in $\Onl$. Consequently, by Lemma \ref{solução trivial local}, the problem 
$$(L_{n\ell}+v_{2}(x)-\lambda_{p}^{\Onl})u=0 \quad \mbox{in }\Onl$$
admits only the trivial solution $u \equiv0$, contradicting the fact that $\varphi_{n\ell}>0$ is a solution to this equation.

Now consider that $\lambda_{p}^{\Onl}$ admits an eigenfunction $\varphi_{n\ell} \in C(\overline{\Omega}_{n\ell}).$ By Lemma \ref{propriedades-NL-autovalor} (iii), we find that $$\lambda_{p}^{\Omega}[\mathfrak{L}+V]=\lambda_{p}^{\Onl}< \min_{\overline{\Omega}_{n\ell}} v_{2}(x).$$  
Therefore, Theorem \ref{Th: caso l_p<m} implies that there exists an eigenfunction $\varphi \in X_{\mathfrak{B}}$ of $\lambda_{p}^{\Omega}[\mathfrak{L}+V]$ such that $\varphi \gg 0$. By the same arguments as before, we have that $\varphi|_{\Onl}$ is a strict positive supersolution for the operator $L_{n\ell}+v_{2}-\lambda_p^{\Omega_{n\ell}}$ in $\Onl$, which again yields a contradiction.
This completes the proof.

\end{proof}

\section{Asymptotic Behavior}\label{sec4}
In this section, we analyze the asymptotic behavior of the local/nonlocal eigenvalue introduced previously. To this end, let 
$$A:=(a_{1},a_{2}) \in L^\infty(\Omega_\ell)\times C(\overline{\Omega}_{n\ell}), \gamma_k:=(\alpha_k,\beta_k) \in C(\overline{\Omega}_\ell)\times C(\overline{\Omega}_{n\ell})$$
and define
\begin{align*}
    \mu_{k}:=\lambda_{p}^{\Omega}[\mathfrak{L}+A+\gamma_{k}]
\end{align*}

The next lemma will be used throughout the subsequent proofs.

\begin{lemma}\label{lema-varphik-local}
    Suppose that $\mu_k= \lambda_p^\Omega[\mathfrak{L}+A+\gamma_k]$ is bounded and admits an eigenfunction $\varphi_k \in X_\mathfrak{B}$ with $\varphi_k\gg0$ and $|\varphi_k|_2=1$. If $\{\alpha_k\}$ is  bounded from below, then there exists a subsequence of $\{\varphi_k|_{\Oe}\}$, still labeled by $k$, such that
    \begin{align*}
\varphi_k \rightharpoonup \varphi_\ell \quad \text{in } H^1_{\Gamma_0}(\Oe) \quad \mbox{and} \quad 
\varphi_k\to \varphi_{\ell}\quad \text{in } L^2(\Oe).
\end{align*}
Furthermore, if there exists an open set $D \subset \Oe$ such that 
\begin{equation}\label{alphak-to-infty}\inf_D \alpha_k \to +\infty,
\end{equation}
then
$$
\varphi_\ell = 0 \quad \mbox{a.e. }x \in D.
$$
\end{lemma}
\begin{proof}
Let $\varphi_k\in X_\mathfrak{B}$, $\varphi_k\gg0$, be the eigenfunction associated with $\mu_k$ such that $|\varphi_k|_2=1$. Multiplying the equation for $\varphi_k|_{\Oe}$ by $\psi \in H_{\Gamma_0}^1(\Oe)$ and integrating by parts over $\Oe$, we obtain   
    \begin{align}
        \int_{\Oe}\nabla \varphi_{k}\nabla \psi+\int_{\Oe}a_{1}(x)\varphi_{k}\psi-\int_{\Oe}\psi \int_{\Onl}J(x-y)\varphi_{k}(y)dydx+ \nonumber \\\int_{\Oe}\alpha_{k}(x)\varphi_{k}\psi=\mu_{k}\int_{\Oe}\varphi_{k}\psi \label{2.9*}
    \end{align}
    Taking $\psi = \varphi_k|_{\Oe}$ and using the boundedness of $\{\alpha_k\}$, we conclude that $\{\varphi_k\}$ is bounded in $H_{\Gamma_0}^1(\Oe)$. Consequently, passing to a subsequence if necessary, we have
    \begin{align*}
\varphi_k \rightharpoonup \varphi_\ell \quad \text{in } H^1_{\Gamma_0}(\Oe) \quad \mbox{and} \quad 
\varphi_k\to \varphi_{\ell}\quad \text{in } L^2(\Oe).
\end{align*}

Now assume that \eqref{alphak-to-infty} holds. Choosing $\psi = \varphi_k|_{\Oe}$ again as a test function in \eqref{2.9*} and using the facts that $\{\mu_k\}$ is bounded, $|\varphi_k|_2 = 1$, and $\{\varphi_k\}$ is bounded in $H_{\Gamma_0}^1(\Oe)$, it follows that there exists a constant $C > 0$ such that 
\begin{equation*}
\int_{\Oe}\alpha_{k}(x)\varphi_{k}^2 \leq C.
\end{equation*}
Let $\underline{\alpha}\in \mathbb{R}$ be such that $\underline{\alpha} \leq \alpha_k(x)$ for all $x \in \Oe$ and $k\geq 1$. Then 
$$
\underline{\alpha}\int_{\Oe \setminus D}\varphi_k^2+ \inf_D\alpha_k\int_D \varphi_k^2 \leq C.
$$
Using \eqref{alphak-to-infty}, we establish that
$$
\int_D \varphi_k^2 \to 0.
$$
By the uniqueness of the limit, $\varphi_\ell = 0$ a.e. in $D$, which completes the proof.
\end{proof}

As a starting point for our asymptotic analysis, we present the proof of Proposition \ref{prop-gammak}.

 \begin{proof}[Proof of Proposition \ref{prop-gammak}]
    First, we consider the case (i). Without loss of generality, we may assume that $\overline{D} \subset \Oe$. Thus, employing the monotonicity properties (Lemmas \ref{propriedades-NL-autovalor}, \ref{propriedades-local-autovalor} and Proposition \ref{estimativas}), we obtain
\begin{align*}
    \mu_{k} < \sigma_{1}^{\Oe}[-\Delta+a_1+\alpha_{k}; \mathfrak{B}]\le \sigma_{1}^{D}[-\Delta+a_1]+\sup_{D}\alpha_{k}(x).
\end{align*}
Therefore, due to \eqref{alphauniK}, the result follows.

Analogously, for case (ii) it is clear that
    \begin{align*}
        \mu_{k}\le \lambda_{p}^{\Onl}[L_{n\ell}+a_{2}+\beta_{k}]\le \lambda_{p}^{D}[L_{D}+a_{2}]+\sup_{D}\beta_{k}(x)
    \end{align*}
    Letting $k \to \infty$ in the above inequality establishes the second result by \eqref{betauniK} and concludes the proof.
 \end{proof}

 \subsection{Proof of Theorem \ref{T1}}
  For ease of reading, we divide the proof of Theorem~\ref{T1} into several propositions.

\begin{proposition}\label{Prop alphk} Let $(\alpha_{k},\beta_{k}) \in C(\overline{\Omega}_{\ell}) \times C(\overline{\Omega}_{n\ell})$ be a non-decreasing sequence. If 
\begin{align} \label{alphauni}
    \alpha_{k}(x) \to + \infty \quad \forall x \in \overline{\Omega}_\ell  \quad \mbox{and} \quad
      \beta_{k} \to \beta^{*} \quad \mbox{in } C(\overline{\Omega}_{n\ell}),
  \end{align}
    then  
    $$ \mu_{k} \to\lambda_{p}^{\Onl}[L_{n\ell}+a_2+\beta^*].$$
    Furthermore, if $\lambda_{p}^{\Onl}[L_{n\ell}+a_2+\beta^*]$ admits a positive principal eigenfunction $\varphi_{n\ell} \in C(\overline{\Omega}_{n\ell})$ normalized in $L^2(\Onl)$, then for sufficiently large $k$, $\mu_k$ admits an eigenfunction $\varphi_{k} \in X_{\mathfrak{B}}$ such that $\varphi_{k} \gg0$, $|\varphi_{k}|_{2}=1$, which satisfies
    \begin{equation}
        \label{conv-eigenfunction-P4.2}\varphi_k \to (0,\varphi_{n\ell}) \quad    \mbox{in } L^2(\Oe)\times C(\overline{\Omega}_{n\ell}).
    \end{equation}

\end{proposition}

\begin{proof}
   Since $(\alpha_{k},\beta_{k})$ is non-decreasing, it follows from Proposition \ref{Propriedades} (ii) that $\{\mu_k\}$ is non-decreasing. Furthermore, 
    by Proposition \ref{estimativas}, we have $$\mu_{k} \le \lambda_{p}^{\Onl}[L_{n\ell}+a_2+\beta_{k}]\le \lambda_{p}^{\Onl}[L_{n\ell}+a_2+\beta^*]=:\lambda_{n\ell}.$$
    Consequently, the following limit exists: 
    \begin{equation}\label{lim-muk}
        \lim_{k \to \infty} \mu_{k}=\mu^*\leq \lambda_{n\ell}.
    \end{equation}
Suppose, by contradiction, that equality does not hold. By the continuity of $\lambda_p^{\Onl}[L_{n\ell}+a_2+\beta_k]$ (see Lemma \ref{propriedades-NL-autovalor}), it follows that
$$\lambda_p^{\Onl}[L_{n\ell}+a_2+\beta_k] \to  \lambda_{n\ell}.$$
Combining this convergence with \eqref{lim-muk} and $\mu^*<\lambda_{n\ell}$, we conclude that for sufficiently large $k\geq1$,
$$\mu_{k}= \lambda_p^{\Omega}[\mathfrak{L}+A(x)+\gamma_k] <  \lambda_p^{\Onl}[L_{n\ell}+a_2+\beta_k]\leq \min_{\overline{\Omega}_{n\ell}}(a_{2}(x)+\beta_k(x)).$$
Therefore, it follows from Theorem \ref{Th: caso l_p<m} that there exists an eigenfunction $\varphi_k\in X_{\mathfrak{B}}$ associated with $\mu_k$ such that $ \varphi_{k} \gg0$ and $|\varphi_{k}|_{2}=1$. By Lemma \ref{lema-varphik-local}, passing to a subsequence if necessary,    \begin{align*}
\varphi_k \rightharpoonup \varphi_\ell \quad \text{in } H^1_{\Gamma_0}(\Oe) \quad \mbox{and} \quad 
\varphi_k\to \varphi_{\ell}\quad \text{in } L^2(\Oe).
\end{align*}
Additionally, since $\alpha_k$ is non-decreasing and \eqref{alphauni} holds, Dini's Theorem yields
$$
\inf_{\Omega_\ell} \alpha_k \to +\infty.
$$
Thus, Lemma \ref{lema-varphik-local} also implies that $\varphi_\ell=0$ a.e. in $\Oe$.

 On the other hand, $\varphi_k|_{\Onl}$ satisfies
\begin{align}\label{2.10}
    L_{n\ell} \varphi_{k}+L_{\ell}\varphi_{k}+(a_2(x)+\beta_{k}(x))\varphi_{k}=\mu_{k}\varphi_{k} \quad \text{in } \Onl.
\end{align}
Since $\mu_{k} \le \mu^{*}<\min_{\overline{\Omega}_{n\ell}} (a_{2}(x)+\beta_{k}(x))\leq \min_{\overline{\Omega}_{n\ell}} (a_{2}(x)+\beta^*(x))$, the equation above is equivalent to
\begin{align}\label{2.11}
    \varphi_{k}=\frac{- L_{n\ell} \varphi_{k}-L_{\ell}\varphi_{k}}{a_{2}(x)+\beta_{k}(x)-\mu_{k}} \quad\text{in } \Onl.
\end{align}
By the continuity of the integral operator, we obtain 
$$- L_{\ell}\varphi_{k}=\int_{\Oe}J(x-y) \varphi_{k}(y)dy \to 0 \quad \mbox{in } C(\overline{\Omega}_{n\ell}).$$
Furthermore, by the compactness of $L_{n\ell}:L^2(\Onl) \to C(\overline{\Omega}_{n\ell})$, for a subsequence still denoted by $k$, we find that
$$-L_{n\ell}\varphi_k=2\int_{\Onl}J(x-y)\varphi_{k}(y)dy \to g \ge 0 \quad \mbox{in } C(\overline{\Omega}_{n\ell}).$$
 As $\mu^*<\lambda_{n\ell} \leq \min(a_2(x)+\beta^*(x))$, letting $k \to \infty$ in \eqref{2.11} yields 
$$ \varphi_{k} \to \frac{g}{a_2(x)+\beta^*(x)-\mu^{*}}=:\varphi_{n\ell} \quad \text{in } C(\overline{\Omega}_{n\ell}).$$
Therefore, passing to the limit as $k\to \infty$ in \eqref{2.10}, we infer that 
\begin{align*}
    L_{n\ell}\varphi_{n\ell}+(a_2(x)+\beta^*(x))\varphi_{n\ell}=\mu^{*}\varphi_{n\ell} \quad \mbox{in }\Onl.
\end{align*}
Noting that 
$$|\varphi_{n\ell}|_{L^{2}(\Omega_{n\ell})}=|(0,\varphi_{n\ell})|_{2}=\displaystyle\lim_{k \to \infty}|\varphi_{k}|_{2}=1,$$
it follows that $\varphi_{n\ell}\ge 0$ and $\varphi_{n\ell}\ne 0$. This implies that  $\varphi_{n\ell}$ is a positive eigenfunction associated with the eigenvalue $\mu^*$. By the uniqueness of the principal eigenvalue, $\mu^*=\lambda_{n\ell}$, which is a contradiction.

Finally, when $\lambda_{p}^{\Onl}[L_{n\ell}+a_2+\beta^*]$ admits a positive eigenfunction, it follows from Lemma \ref{propriedades-NL-autovalor} (iii) that 
$$
\lambda_{p}^{\Onl}[L_{n\ell}+a_2+\beta^*] < \min_{\overline{\Omega}_{n\ell}}(a_2(x)+\beta^*(x)).
$$
Moreover,
$$
\mu_k \to \lambda_{p}^{\Onl}[L_{n\ell}+a_2+\beta^*] \quad \mbox{and} \quad m_k:=\min_{\overline{\Omega}_{n\ell}}(a_2(x)+\beta_k(x))\to \min_{\overline{\Omega}_{n\ell}}(a_2(x)+\beta^{*}(x)).
$$
Thus, for sufficiently large $k\geq 1$, we have $\mu_k<m_k$. By Theorem \ref{Th: caso l_p<m}, $\mu_k$ admits a principal eigenfunction, and the arguments above yield the convergence of the eigenfunctions. Since this reasoning applies to any subsequence, \eqref{conv-eigenfunction-P4.2} holds. This concludes the proof.
  
\end{proof}
 Next we consider the case where $(\alpha_{k},\beta_{k}) \to (\alpha^*,+\infty)$.

\begin{proposition}\label{lambdk infito}
Let  $(\alpha_{k},\beta_{k}) \in C(\overline{\Omega}_{\ell}) \times C(\overline{\Omega}_{n\ell})$ be a non-decreasing sequence. If
\begin{align}\label{betauniO}    \alpha_{k} \to \alpha^* \quad \mbox{in } C(\overline{\Omega}_{\ell}) \quad \mbox{and} \quad    \beta_{k}(x) \to + \infty \quad \forall x\in\overline{\Omega}_{n\ell},
    \end{align}
   then
   $$ \mu_{k} \to \sigma_{1}^{\Oe}[-\Delta+a_1+\alpha^*;\mathfrak{B}].$$
    Furthermore, for sufficiently large $k$, $\mu_k$ admits an eigenfunction $\varphi_{k} \in X_{\mathfrak{B}}$ such that $\varphi_{k} \gg0$, $|\varphi_{k}|_{2}=1$, which satisfies
   $$\varphi_k \to (\varphi_\ell,0) \quad\mbox{in } X_\mathfrak{B}=W^{2,p}_\mathfrak{B}(\Oe) \times C(\overline{\Omega}_{n\ell}),$$
    where $\varphi_\ell\gg0$ denotes the positive eigenfunction associated with $\sigma_1^{\Omega_\ell}[-\Delta +a_1 +\alpha^*, \mathfrak{B}]$, normalized such that $|\varphi_\ell|_{L^2(\Oe)}=1$.
\end{proposition}

\begin{proof}

    Since $(\alpha_k, \beta_k)$ is non-decreasing, it follows from Proposition \ref{Propriedades} (ii) that $\{\mu_k\}$ is monotonically non-decreasing. Additionally, 
    by the monotonicity properties (see Lemma \ref{propriedades-NL-autovalor} and Proposition \ref{estimativas}), we have
$$\mu_{k}< \sigma_{1}^{\Oe}[-\Delta+a_{1}+\alpha_{k};\mathfrak{B}] \leq\sigma_{1}^{\Oe}[-\Delta+a_{1}+\alpha^*;\mathfrak{B}]=: \sigma_1.$$
    Therefore, the following limit exists:
    $$\lim_{k\to \infty}\mu_{k}=\mu^{*}\le \sigma_{1}.$$
    On the other hand, it follows from \eqref{betauniO} that $$\displaystyle \min_{\overline{\Omega}_{n\ell}}(a_2(x)+\beta_{k}(x)) \to + \infty.$$
    Since $\{\mu_k\}$ is bounded, for sufficiently large $k\geq1$, we have
    \begin{equation}\label{muk cota}
        \mu_{k}<\displaystyle \min_{\overline{\Omega}_{n\ell}}(a_2(x)+\beta_{k}(x)).
    \end{equation}
    Thus, Theorem \ref{Th: caso l_p<m} guarantees the existence of $\varphi_{k} \in X_{\mathfrak{B}}$ such that $\varphi_{k} \gg 0$ and $|\varphi_{k}|_{2}=1$. By Lemma \ref{lema-varphik-local}, passing to a subsequence if necessary,       \begin{align*}
\varphi_k \rightharpoonup \varphi_\ell \quad \text{in } H^1_{\Gamma_0}(\Oe) \quad \mbox{and} \quad 
\varphi_k\to \varphi_{\ell}\quad \text{in } L^2(\Oe).
\end{align*}
Furthermore, $\varphi_{k}|_{\Onl}$ also satisfies 
\begin{align*}
    L_{n\ell}\varphi_{k}+L_{\ell}\varphi_{k}+(a_{2}(x)+\beta_{k}(x))\varphi_{k}=\mu_{k}\varphi_{k}  \quad \mbox{in } \Onl.
\end{align*}
In view of \eqref{muk cota}, we can rewrite this equality as
\begin{equation}
\label{varphik-parte-nl}\varphi_{k}=\frac{ -L_{n\ell}\varphi_{k}-L_{\ell}\varphi_{k}}{a_2(x)+\beta_{k}(x)-\mu_{k}} \quad \mbox{in }\Onl.
\end{equation}
As $|\varphi_k|_2=1$, the properties of the integral operator imply that there exists $C>0$ such that
$$\|-L_{n\ell}\varphi_{k}-L_{\ell}\varphi_{k}\|_{C(\overline{\Omega}_{n\ell})}\leq C \quad \forall k.$$
Letting $k \to +\infty$ in \eqref{varphik-parte-nl} and using \eqref{betauniO}, we find that
$$
\varphi_k \to 0 \quad \mbox{in }C(\overline{\Omega}_{n\ell}).
$$

On the other hand, letting $k \to \infty$ in the weak formulation for $\varphi_k|_{\Oe}$ (see \eqref{2.9*}) leads to
    \begin{align*}
        \int_{\Oe}\nabla \varphi_{\ell}\nabla \psi+\int_{\Oe}a_{1}(x)\varphi_{\ell}\psi+ \int_{\Oe}\alpha^*(x)\varphi_{\ell}\psi=\mu^{*}\int_{\Oe}\varphi_{\ell}\psi \quad \forall\psi \in H_{\Gamma_0}^1(\Oe).
    \end{align*}
Noting that
$$|\varphi_{\ell}|_{L^{2}(\Omega_{\ell})}=|(\varphi_\ell,0)|_{2}=\displaystyle\lim_{k \to \infty}|\varphi_{k}|_{2}=1,$$ we have $\varphi_{\ell} \ge0$ and $\varphi_{\ell} \ne 0$. Thus, $\varphi_\ell$ is a positive eigenfunction associated with the eigenvalue $\mu^*$.  
    By the uniqueness of the principal eigenvalue, we conclude that $\mu^{*}=\sigma_{1}$, which contradicts the assumption that $\mu^*<\sigma_1$. 
     
     To complete the proof, it remains to show that $$\varphi_{k} \to \varphi_{\ell} \quad \mbox{in } W_\mathfrak{B}^{2,p}(\Oe).$$  Indeed, we know that $\varphi_{k}|_{\Oe}$ satisfies
    \[
\begin{cases}
-\Delta \varphi_{k}+ a_1(x)\varphi_{k}+\alpha_{k}(x)\varphi_{k}-\mu_{k} \varphi_{k}=v_k:=  \int_{\Onl}J(x-y)\varphi_{k}(y) \,dy& \text{in } \Omega_{\ell}, \\ 
\mathfrak{B}\varphi_{k} = 0 & \text{on } \partial \Omega_{\ell}. \\ 

\end{cases} 
\]
Since $\{v_k\}$ is a bounded sequence in $C(\overline{\Omega}_{\ell})$, elliptic regularity results imply that $\{\varphi_k\}$ is bounded in $ W^{2,p}_{\mathfrak{B}}(\Oe)$. Compact embeddings ensure that, passing to a subsequence if necessary, $\varphi_k \to \varphi_\ell$ in $C^{1,\theta}(\overline{\Omega}_\ell)$ for some $0<\theta<1$. Finally, noting that
    \[
\begin{cases}
-\Delta (\varphi_{k}-\varphi_\ell)+ a_1(x)(\varphi_{k}-\varphi_\ell)= w_k & \text{in } \Omega_{\ell}, \\ 
\mathfrak{B}(\varphi_{k} -\varphi_\ell)= 0 & \text{on } \partial \Omega_{\ell}. \\ 

\end{cases} 
\]
where 
$$
w_k := (\mu_k-\alpha_{k}(x))(\varphi_{k}- \varphi_\ell)+ (\mu_{k}-\mu^{*})\varphi_{\ell}-(\alpha_{k}(x)-\alpha^{*}(x))\varphi_{\ell}+ \int_{\Onl}J(x-y)\varphi_{k}(y) \,dy.
$$
Since $w_k \to 0$ in $L^p(\Oe)$ for all $p>1$, applying an elliptic regularity result once more concludes the proof.

\end{proof}

Finally, we consider the scenario where both parameters diverge to infinity simultaneously.

\begin{proposition}\label{T1(iii)}
   Let $(\alpha_{k},\beta_{k}) \in C(\overline{\Omega}_{\ell}) \times C(\overline{\Omega}_{n\ell})$ be a non-decreasing sequence.
   If \begin{align}
           \alpha_{k}(x) &\to +\infty \quad \forall x \in  \overline{\Omega}_{\ell} \label{alphak1}  \quad\mbox{and}\\ \beta_{k}(x) &\to +\infty \quad \forall x \in \overline{\Omega}_{n\ell},\label{betak1}
       \end{align}
       then
    $$ \mu_{k} \to+\infty.$$
\end{proposition}
\begin{proof}
Arguing as in the previous proofs, $\{\mu_k\}$ is non-decreasing and the limit exists:
    $$
    \lim_{k \to \infty} \mu_{k}=\mu^{*} \leq+\infty.
    $$
    Suppose, by contradiction, that $\mu^*<\infty$.
   As in Proposition \ref{lambdk infito}, we have $$\mu_{k} = \lambda_p^{\Omega}[\mathfrak{L}+A+\gamma_k]\leq\displaystyle \min_{\overline{\Omega}_{n\ell}}(a_{2}(x)+\beta_{k}(x))\to +\infty.$$
    Since $\mu_k\leq \mu^*$, for sufficiently large $k\geq1$, we have $$\mu_{k}<\displaystyle \min_{\overline{\Omega}_{n\ell}}(a_{2}(x)+\beta_{k}(x)) .$$
    By Theorem \ref{Th: caso l_p<m}, we guarantee the existence of an eigenfunction $\varphi_{k} \in X_{\mathfrak{B}} $ such that $\varphi_{k} \gg0$ and $|\varphi_{k}|_{2}=1$. Applying Lemma \ref{lema-varphik-local}, we conclude that, up to a subsequence,       \begin{align*}
\varphi_k \rightharpoonup \varphi_\ell \quad \text{in } H^1_{\Gamma_0}(\Oe) \quad \mbox{and} \quad 
\varphi_k\to \varphi_{\ell}\quad \text{in } L^2(\Oe).
\end{align*}
Since $\{\alpha_k\}$ is non-decreasing and \eqref{alphak1} holds, Dini's Theorem yields
$$
\inf_{\Omega_\ell} \alpha_k \to +\infty.
$$
Thus, Lemma \ref{lema-varphik-local} also implies that $\varphi_\ell=0$ a.e. in $\Oe$.

 Furthermore, since \eqref{betak1} holds, we can repeat the arguments from the proof of Proposition \ref{lambdk infito} to conclude that 
    $$\varphi_{k} \to 0 \quad \mbox{in  } C(\overline{\Omega}_{n\ell}).$$ 
    Consequently,
    $$
    \varphi_k \to 0 \quad \mbox{in }L^2(\Omega),
    $$
    which contradicts the fact that $|\varphi_k|_2=1$. This concludes the proof.

\end{proof}

The proof of Theorem \ref{T1} then follows directly from Propositions \ref{Prop alphk}, \ref{lambdk infito}, and \ref{T1(iii)}.

\subsection{Proof of Theorem \ref{T2}} \label{sec4.2}

We now turn our attention to the proof of Theorem \ref{T2}, which deals with the case where the  weights  vanish. As in the previous subsection, we divide the proof into three propositions. The first result establishes the asymptotic behavior of $\mu_{k}$ when $\alpha_{k}=0$ in a subdomain $\Omega_0 \subset \Oe$. More precisely:

\begin{proposition}\label{Refugio em local}
        Let $(\alpha_{k},\beta_{k}) \in C(\overline{\Omega}_{\ell}) \times C(\overline{\Omega}_{n\ell})$ be a non-decreasing sequence. Suppose there exists a subdomain $\Omega_0$,  $\overline{\Omega}_0 \subset \Oe$, such that
       \begin{equation}
            \alpha_k(x) = 0 \quad \forall x \in \Omega_0, k\geq1 \quad \mbox{and} \quad\alpha_{k}(x)\to+\infty \quad \forall x \in \Omega_{\ell}\setminus\overline{\Omega}_0, \label{alpha-to-infty1}
        \end{equation}
        and that 
         \begin{equation}\label{beta-to-infinity}
            \beta_{k}(x) \to +\infty  \quad \forall x \in \overline{\Omega}_{n\ell}.
        \end{equation}
        Then $$\mu_{k} \to \sigma_{1}^{\Omega_{0}}[-\Delta+a_{1}].$$
        Furthermore, for sufficiently large $k\geq1$, $\mu_k$ admits an eigenfunction $\varphi_{k} \in X_{\mathfrak{B}} $ such that $\varphi_{k} \gg0$, $|\varphi_{k}|_{2}=1$, satisfying
        $$\varphi_{k} \to \phi \quad \mbox{in } H_{\Gamma_{0}}^{1}(\Oe)$$ where
        $$ \phi:= \left\{ 
        \begin{array}{ll}
            \varphi_0 & \mbox{in }\Omega_0,  \\
             0& \mbox{in }\Omega_{\ell}\setminus \overline{\Omega}_0,
        \end{array}
        \right.$$
        and $\varphi_0 \in W^{2,p}_0(\Omega_0)$, for $p>N$, denotes the positive eigenfunction associated with $\sigma_1^{\Omega_0}[-\Delta+a_1]$, normalized such that $|\varphi_0|_{L^2(\Omega_0)}=1$. Moreover,
        $$
        \varphi_k \to \varphi_0 \quad \mbox{in }L^\infty_{loc}(\Omega_0).
        $$
\end{proposition}
\begin{proof}
    Applying the monotonicity properties of eigenvalues with respect to the domain, we conclude that
    $$\mu_{k} < \sigma_{1}^{\Oe}[-\Delta+a_{1}+\alpha_{k}(x);\mathfrak{B}] < \sigma_1^{\Omega_0}[-\Delta+a_1]=:\sigma_0.$$
    Additionally, because $(\alpha_k,\beta_k)$ is non-decreasing, $\{\mu_k\}$ is monotonically non-decreasing, and therefore the limit exists:
    \begin{align*}
        \lim_{k\to \infty}\mu_{k}=\mu^{*} \le \sigma_{0}
    \end{align*}
    Let us prove that equality holds. Using  \eqref{beta-to-infinity}, it follows that $$\displaystyle \min_{\overline{\Omega}_{n\ell}}(a_2(x)+\beta_{k}(x))  \to +\infty.$$
    Since $\{\mu_{k}\}$ is bounded, the convergence above implies that for sufficiently large $k \ge 1$,
    \begin{align}\label{estimina muk 2}
        \mu_{k} < m_{k}:=\displaystyle \min_{\overline{\Omega}_{n\ell}}(a_2(x)+\beta_{k}(x)).
    \end{align}
    Thus, Theorem \ref{Th: caso l_p<m} guarantees that $\mu_k$ admits an eigenfunction $\varphi_{k} \in X_{\mathfrak{B}} $ such that $\varphi_{k} \gg0$ and $|\varphi_{k}|_{2}=1$. By Lemma \ref{lema-varphik-local}, passing to a subsequence if necessary,       \begin{align*}
\varphi_k \rightharpoonup \varphi_\ell \quad \text{in } H^1_{\Gamma_0}(\Oe) \quad \mbox{and} \quad 
\varphi_k\to \varphi_{\ell}\quad \text{in } L^2(\Oe).
\end{align*}
Let $B \subset \Oe \setminus \overline{\Omega}_{0}$ be an arbitrary open ball such that $\overline{B} \subset \Oe \setminus \overline{\Omega}_{0}$. Because $\alpha_k$ is non-decreasing and \eqref{alpha-to-infty1} holds, we have
$$
\inf_{B} \alpha_k \to +\infty.
$$
Applying Lemma \ref{lema-varphik-local} again, we conclude that $\varphi_{\ell}=0$ a.e. in $B$ and, consequently, $\varphi_\ell=0$ a.e. in $\Oe \setminus\overline{ \Omega}_{0}$.  Since $\Omega_0$ is a regular subdomain with $\overline{\Omega}_0 \subset \Oe$, we have $\varphi_{\ell} \in H_{0}^{1}(\Omega_{0})$ (see, for instance, \cite[Theorem 9.17]{brezis2010functional}). On the other hand, $\varphi_{k}|_{\Onl}$ satisfies 
$$L_{n\ell}\varphi_{k}+L_{\ell}\varphi_{k}+a_2(x)\varphi_{k}+\beta_{k}(x)\varphi_{k}=\mu_{k}\varphi_{k}  \quad \mbox{in } \Onl.$$
Arguing as in the proof of Proposition \ref{lambdk infito}, we ensure that
$$\varphi_{k} \to 0 \ \quad \mbox{in } C(\overline{\Omega}_{n\ell}).$$
Finally, since that $\alpha_k=0$ in $\Omega_{0}$,  multiplying the equation for $\varphi_k|_{\Oe}$ by $\psi \in H_{0}^1(\Omega_{0})$ and integrating by parts over $\Omega_0$, we obtain
\begin{align}\label{solfraca-P2.16}
     \int_{\Omega_0}\nabla \varphi_{k}\nabla \psi+\int_{\Omega_{0}}a_{1}(x)\varphi_{k}\psi-\int_{\Omega_0}\psi \int_{\Onl}J(x-y)\varphi_{k}(y)dydx=\mu_{k}\int_{\Omega_0}\varphi_{k}\psi.
\end{align}
Using the established convergences and letting $k \to \infty$ in the identity above, we conclude that 
$$\int_{\Omega_0}\nabla \varphi_{\ell}\nabla \psi+\int_{\Omega_{0}}a_{1}(x)\varphi_{\ell}\psi=\mu^*\int_{\Omega_0}\varphi_{\ell}\psi \quad  \forall \psi \in H_0^1(\Omega_0).
$$
Since
$$|\varphi_{\ell}|_{L^{2}(\Omega_{\ell})}=|(\varphi_\ell,0)|_{2}=\displaystyle\lim_{k \to \infty}|\varphi_{k}|_{2}=1,$$ we have  $\varphi_{\ell} \ge0$ and $\varphi_{\ell} \ne 0$. By the uniqueness of the principal eigenvalue, we find that $\mu^*=\sigma_{0}$ and $\varphi_{\ell} =\varphi_{0}$. Now, we prove that $$\varphi_{k} \to \phi \quad \mbox{in } H^{1}(\Oe). $$
By multiplying the $\varphi_{k}|_{\Oe}$ equation by $\psi=\varphi_{k}-\varphi_{\ell}$ and integrating by parts in $\Oe$, we obtain  
\begin{align}\label{2.333}
    \int_{\Oe}\nabla \varphi_{k}\nabla(\varphi_{k}-\varphi_{\ell}) +\int_{\Oe}a_{1}(x)\varphi_{k}(\varphi_{k}-\varphi_{\ell})+\int_{\Oe}\alpha_{k}(x)\varphi_{k}(\varphi_{k}-\varphi_{\ell})-\\\int_{\Oe}\int_{\Onl}J(x-y)\varphi_{k}(y)(\varphi_{k}-\varphi_{\ell})dydx \nonumber =\mu_{k}\int_{\Oe}\varphi_{k}(\varphi_{k}-\varphi_{\ell})
\end{align}
Using the previous limits, we conclude that
\begin{align*}
    \int_{\Oe}a_{1}(x)\varphi_{k}(\varphi_{k}-\varphi_{\ell}) \to 0, \quad \int_{\Oe}\int_{\Onl}J(x-y)\varphi_{k}(y)(\varphi_{k}-\varphi_{\ell})dydx \to 0
\end{align*}
and
\begin{align*}
    \mu_{k}\int_{\Oe}\varphi_{k}(\varphi_{k}-\varphi_{\ell}) \to 0.
\end{align*}
Consequently, since $\alpha_{k}= 0 \ \mbox{for all }  x \in \Omega_{0}$, $k \ge 1$ and $\varphi_{\ell}=0 $ a.e. in $\Oe \setminus \overline{\Omega}_{0}$, \eqref{2.333} can be rewritten as
\begin{align*}
    \int_{\Oe}|\nabla(\varphi_{k}-\varphi_{\ell})|^{2}+\int_{\Oe}\nabla \varphi_{\ell}\nabla(\varphi_{k}-\varphi_{\ell})+\int_{\Oe \setminus\Omega_{0}}\alpha_k(x)\varphi_{k}^{2}\leq\operatorname{o}_{k}(1).
\end{align*}
Letting $k \to +\infty$ in this inequality, we obtain that $\varphi_k \to \varphi_\ell$ in $H^1(\Oe)$.

It remains to show that 
$$\varphi_{k} \to \varphi_{0} \quad\mbox{in } L^\infty_{loc}(\Omega_{0}).$$
Indeed, in $\Omega_{0}$ we have 
\begin{multline*}
      -\Delta (\varphi_{k}-\varphi_{0}) = (\mu_{k}-\sigma_{1}^{\Omega_{0}})\varphi_{k} + \sigma_{1}^{\Omega_{0}}(\varphi_k-\varphi_{0}) \\
- a_{1}(x)(\varphi_{k}-\varphi_{0}) + \int_{\Onl}J(x-y)\varphi_{k}(y)dy =: f_{k} \quad \mbox{in } \Omega_{0}.
\end{multline*}
From the limits obtained previously, we have that
$$f_k \to 0 \quad\mbox{in } L^{2^*}(\Omega_{0}).$$
If $N \leq 2$, the above convergence holds for all $1 < p < +\infty$. By the Agmon-Douglis-Nirenberg interior estimates (\cite[Theorem 9.11]{gilbarg1998elliptic}), it follows that for any subdomain $\omega_{1} \subset\subset \Omega_{0}$, there exists a constant $C>0$ such that
$$
\|\varphi_k-\varphi_0\|_{W^{2,p}(\omega_1)} \leq C\left(\|\varphi_{k}-\varphi_{0}\|_{L^{p}(\Omega_0)}+\|f_{k}\|_{L^{p}(\Omega_0)}\right),
$$
therefore $\varphi_k \to \varphi_0$ in $W^{2,p}_{loc}(\Omega_0)$, with $1 < p < +\infty$. By the embedding $W^{2,p}_{loc}(\Omega_0) \hookrightarrow C^{0,\theta}_{loc}(\Omega_0)$, for $p > N/2$ and some $0 < \theta < 1$, we obtain the result.  

If $N \geq 3$, the preceding argument guarantees that $\varphi_k \to \varphi_0$ in $W^{2,2^*}_{loc}(\Omega_0) \hookrightarrow L^{p_1}_{loc}(\Omega_0)$, where $p_1$ is given by $1/p_1 = 1/2^* - 2/N$. If $p_1 > N/2$, then
$$
f_k \to 0 \quad \mbox{in } L^{p_1}_{loc}(\Omega_0)
$$
and again, interior estimates allow us to conclude that $\varphi_k \to \varphi_0$ in $W^{2,p_1}_{loc}(\Omega_0) \hookrightarrow C^{0,\theta}_{loc}(\Omega_0)$, for some $0 < \theta < 1$, concluding the proof. In case $p_1 < N/2$, we can repeat this bootstrap argument a finite number of times to conclude that $\varphi_k \to \varphi_0$ in $W^{2,p_m}_{loc}(\Omega_0)$, where $1/p_m = 1/2^* - 2m/N$ with $m \in \mathbb{N}$ such that $p_m > N/2$. This completes the proof.
\end{proof}

Next, we establish the asymptotic behavior when $\beta_{k}=0$ in some subdomain of $\Onl$. To this end, given a subdomain $\Omega_0 \subset\Onl$, we define the operator $L_0:C(\overline{\Omega}_0) \to C(\overline{\Omega}_0)$ by
$$L_{0}u(x):=-2\int_{\Omega_{0}}J(x-y)u(y)\, dy.$$

\begin{proposition}\label{refugio nao local}
Let $(\alpha_{k},\beta_{k})$ be a non-decreasing sequence. Suppose there exists a subdomain $\Omega_0 \subset \Onl$ such that
         \begin{equation}
        \beta_k(x) = 0 \quad \forall x \in \overline{\Omega}_0, k\geq1 \quad \mbox{and} \quad     \beta_{k}(x) \to +\infty \ \ \quad \forall x \in \overline{\Omega}_{n\ell}\setminus\overline{\Omega}_0,\label{beta-to-infinto 2}
        \end{equation}
         
        and  \begin{equation}\label{alpha-to-infinity}
            \alpha_{k}(x) \to +\infty  \quad \forall x \in  \overline{\Omega}_{\ell},
        \end{equation}
then $$\mu_{k} \to \lambda_{p}^{\Omega_{0}}[L_0+a_2].$$
    Furthermore, if $\lambda_{p}^{\Omega_0}[L_{0}+a_2]$ admits a positive principal eigenfunction $\varphi_0 \in C(\overline{\Omega}_0)$ normalized in $L^2(\Omega_0)$, then for sufficiently large $k$, $\mu_k$ admits an eigenfunction $\varphi_{k} \in X_{\mathfrak{B}}$ such that $\varphi_{k} \gg0$, $|\varphi_{k}|_{2}=1$, satisfying
    $$\varphi_k \to (0,\phi) \quad    \mbox{in } L^2(\Oe)\times L^2(\Onl),$$
    where 
$$\phi:= \left\{ 
        \begin{array}{ll}
            \varphi_0 & \mbox{in }{\Omega}_0,  \\
             0& \mbox{in }\Omega_{n\ell}\setminus \Omega_0.
        \end{array}
        \right.$$
        Moreover,
        $$
        \varphi_k|_{\overline{\Omega}_0} \to \varphi_0 \quad \mbox{in }C(\overline{\Omega}_0).
        $$
\end{proposition}
\begin{proof}
    Applying the monotonicity properties with respect to the domain, we conclude that
    $$\mu_{k} \le \lambda_{p}^{\Onl}[L_{n\ell}+a_2+\beta_{k}]\le \lambda_{p}^{\Omega_{0}}[L_0+a_2]=: \lambda_{0}.$$
    Because $(\alpha_{k},\beta_{k})$ is non-decreasing, $\mu_k$ is also non-decreasing, and the following limit exists 
    $$\lim_{k \to + \infty} \mu_{k}= \mu^{*} \le \lambda_{0} \le \min_{\overline{\Omega}_{0}} a_2(x).$$ Suppose, by contradiction, that $\mu^{*} < \lambda_{0}$.      The condition \eqref{beta-to-infinto 2} ensures that  \begin{align}\label{convergencia minim}
        \displaystyle\min_{\overline{\Omega}_{n\ell}} (a_2(x)+\beta_{k}(x)) \to \min_{\overline{\Omega}_{0}} a_2(x) > \mu^*
    \end{align}
    Thus, for sufficiently large $k\ge1$, 
    \begin{align}\label{estimitiva mu*}
            \mu_{k} \le \mu^{*}<\min_{\overline{\Omega}_{n\ell}} (a_2(x)+\beta_{k}(x)) 
    \end{align}
   Consequently, by Theorem \ref{Th: caso l_p<m}, there exists an eigenfunction $\varphi_{k} \in X_{\mathfrak{B}}$ associated with $\mu_{k}$ such that $\varphi_{k} \gg 0$ and $|\varphi_{k}|_{2}=1$.  Applying Lemma \ref{lema-varphik-local}, we find that, up to a subsequence,       \begin{align*}
\varphi_k \rightharpoonup \varphi_\ell \quad \text{in } H^1_{\Gamma_0}(\Oe) \quad \mbox{and} \quad 
\varphi_k\to \varphi_{\ell}\quad \text{in } L^2(\Oe).
\end{align*}
Since $\alpha_k$ is non-decreasing and \eqref{alpha-to-infinity} holds, Dini's Theorem yields
$$
\inf_{\Omega_\ell} \alpha_k \to +\infty.
$$
Thus, Lemma \ref{lema-varphik-local} also implies that $\varphi_\ell=0$ a.e. in $\Oe$.

On the other hand, in view of \eqref{estimitiva mu*}, $\varphi_{k}|_{\Onl}$  satisfies
\begin{align} \varphi_{k}=\frac{-L_{n\ell}\varphi_{k}-L_{\ell}\varphi_{k}}{a_2(x)+\beta_{k}(x)-\mu_{k}} \quad \mbox{in } \Onl.\label{2.30}
\end{align}
Arguing analogously to Proposition \ref{Refugio em local}, we obtain 
\begin{align*}
   - L_{n\ell}\varphi_{k} \to g \ge 0  \quad \text{in } C(\overline{\Omega}_{n\ell}) \quad \mbox{and } L_{\ell}\varphi_{k} \to 0 \quad \mbox{in } C(\overline{\Omega}_{n\ell}).
\end{align*}
 In particular, there exists $M>0$ (independent of $k$) such that 
 \begin{equation}\label{limitacao-superior-numerador}
     \|-L_{n\ell}\varphi_{k}-L_{\ell}\varphi_{k}\|_{C(\overline{\Omega}_{n\ell})}\le M.
 \end{equation}
Consequently, letting $k \to +\infty$ in \eqref{2.30} and using \eqref{beta-to-infinto 2}, we deduce that
\begin{align}
    \varphi_{k} \to 0 \quad \mbox{a.e. in } \Onl \setminus \overline{\Omega}_{0}.\label{2.32}
\end{align}
Furthermore, because
$$
a_2(x) + \beta_k(x) - \mu_k \geq \min_{\overline{\Omega}_{n\ell}}(a_2(x)+\beta_k(x))-\mu^* \to \min_{\overline{\Omega
}_0}a_2(x)-\mu^*>0,
$$
by fixing $0<\varepsilon_0< \min_{\overline{\Omega}_0}a_2(x)-\mu^*$, we have for sufficiently large $k\geq1$ that
\begin{equation}\label{limitacao-inferior-denominador}
    \varepsilon_0\leq  a_2(x) + \beta_k(x) - \mu_k \quad \forall x \in \overline{\Omega}_{n\ell}.
\end{equation}
Using \eqref{limitacao-superior-numerador} and \eqref{limitacao-inferior-denominador} in \eqref{2.30}, we conclude that
\begin{align}
    \varphi_{k} \le \frac{M}{\varepsilon_{0}} \in L^{2}(\Onl \setminus \Omega_{0}). \label{2.33}
\end{align}
By   \eqref{2.32} and \eqref{2.33}, we can apply the Dominated Convergence Theorem to infer that 
$$ \varphi_{k} \to 0 \quad \mbox{in } L^{2}(\Onl \setminus \Omega_{0}).$$
Subsequently,  since $\beta_k(x)=0$ for all $x \in \Omega_0$ and $k\geq1$, we have 
\begin{align}
\varphi_{k}=\frac{-L_{n\ell}\varphi_{k}-L_{\ell}\varphi_{k}}{a_{2}(x)-\mu_{k}} \quad \mbox{in } \overline{\Omega}_{0}. \label{Omega0nl}
\end{align}
Letting $k \to + \infty$ in the identity above makes it apparent that 
$$  \varphi_{k}  \to \frac{g}{a_2(x)-\mu^*}=: \varphi_{n\ell} \quad \mbox{in } C(\overline{\Omega}_{0}).$$
By the continuity of the integral operator, we obtain that
\begin{multline*}
    L_{n\ell} \varphi_k =- 2\int_{\Onl\setminus \Omega_0} J(x-y) \varphi_k(y) \, dy - 2\int_{\Omega_0} J(x-y) \varphi_k(y) \, dy \\
    \to -2\int_{\Omega_0} J(x-y) \varphi_{n\ell}(y) \, dy=L_0 \varphi_{n\ell} \quad \mbox{in }C(\overline{\Omega}_0).
\end{multline*}
Thus, letting $k \to +\infty$ in \eqref{Omega0nl}, we find that
\begin{align*}
    \varphi_{n\ell}= \frac{-L_0\varphi_{n\ell}}{a_{2}(x)-\mu^{*}} \quad \mbox{in }\Omega_0.
\end{align*}
Since 
$$  1=\lim_{k \to + \infty}|\varphi_{k}|_{2}^{2}=\lim_{k \to +\infty}\left(\int_{\Omega_0}\varphi_{k}^{2}+\int_{\Onl \setminus \Omega_{0}}\varphi_{k}^{2}+\int_{\Oe}\varphi_{k}^{2}\right)=|\varphi_{n\ell}|^2_{L^2(\Omega_0)},$$
we have  $\varphi_{n\ell} \ge0$ and $\varphi_{n\ell} \ne 0$.  By the uniqueness of the principal eigenvalue, we obtain $\mu^{*}=\lambda_{0}$ and $\varphi_{n\ell}=\varphi_0$, which contradicts the assumption that $\mu^{*}<\lambda_{0}$. 

Finally, assuming that $\lambda_p^{\Omega_0}[L_0+a_2]$ admits a positive eigenfunction, Lemma \ref{propriedades-NL-autovalor} (iii) implies that 
$$
\lambda_{p}^{\Omega_0}[L_{0}+a_2] < \min_{\overline{\Omega}_0}a_2.
$$
Moreover,
$$
\mu_k \to \lambda_{p}^{\Omega_0}[L_{0}+a_2] \quad \mbox{and} \quad m_k:=\min_{\overline{\Omega}_{n\ell}}(a_2+\beta_k)\to \min_{\overline{\Omega}_0}a_2.
$$
Thus, for sufficiently large $k\geq1$, we have $\mu_k<m_k$. By Theorem \ref{Th: caso l_p<m}, $\mu_k$ admits a principal eigenfunction, and we can repeat the above arguments to obtain the result. This concludes the proof.

\end{proof}

Finally, we consider the case where both $\alpha_k$ and $\beta_k$ vanish in some subset.
 To this end,  let  $\Omega_{\ell}^0 \subset \Oe$ and $\Omega_{n\ell}^0 \subset \Onl$ be regular subdomains such that $\overline{\Omega}_\ell^0 \subset\Omega_\ell$.  Setting $\Omega_0=\Omega_{\ell}^0 \cup \Omega_{n\ell}^0 \subset \Omega$, we define the sets
\[
    X_0= \{u \in L^2(\Omega_0) : u|_{\Omega_{\ell}^{0}} \in W_0^{2,p}(\Omega_{\ell}^0),\, p > N, \text{ and }  u|_{\overline{\Omega}_{n\ell}^0} \in C(\overline{\Omega}_{n\ell}^0) \},
\]
$$
Y_0=\{u \in L^2(\Omega_0); \, u|_{\Oe^{0}} \in L^p(\Oe^0), p>N, \text{ and }  u|_{\overline{\Omega}_{n\ell}^0}  \in C(\overline{\Omega}_{n\ell}^0)\}
$$

 For the next result, we employ the following notation: the operators $L_{n\ell}^{0}:C(\overline{\Omega}_{n\ell}^{0}) \to C(\overline{\Omega}_{n\ell}^{0})$ and $L_{\ell}^{0}:L^p({\Omega}_{\ell}^{0}) \to C(\overline{\Omega}_{n\ell}^{0})$ are given by

$$L_{n\ell}^{0}u :=-2\int_{\Omega_{n\ell}^{0}}J(x-y) u(y)\, dy \quad \mbox{and} \quad L_{\ell}^{0}u :=-\int_{\Omega_{\ell}^{0}}J(x-y)u(y)\, dy.$$
Furthermore,
\begin{align*}
    \mathcal{L}_1^{0}u(x)& = -\Delta u - \int_{\Omega_{n\ell}^{0}} J(x-y)u(y)dy  \\ \mathcal{L}_2^{0}u(x) &= -2\int_{\Omega_{n\ell}^0} J(x-y)u(y)dy - \int_{\Oe^0} J(x-y)u(y)dy=L_{n\ell}^{0}u+L_{\ell}^{0}u.
\end{align*} 
Finally, the operator $\mathfrak{L}_{0}:X_0 \to Y_0$ is given by
\[
 \mathfrak{L}_{0}u:= \begin{cases} \mathcal{L}_1^{0} u  & \text{in } \Omega_{\ell}^{0}, \\ \mathcal{L}_2^{0} u  & \text{in } \Omega_{n \ell}^{0}. \end{cases}
\]
Figure \ref{fig:2} illustrates two admissible configurations for the subdomains $\Oe^{0}$ and $\Onl^{0}$.

\begin{proposition}\label{refugio local e nlocal}
Let $(\alpha_{k},\beta_{k})$ be a non-decreasing sequence and let $\Omega_{\ell}^0 \subset \Oe$ and $\Omega_{n\ell}^0 \subset \Onl$ be  as above. Suppose that
        $$
        \alpha_k(x) = 0 \quad \forall x\in \Oe^0 \quad \mbox{and} \quad \beta_k(x) = 0 \quad  \forall x \in \Onl^0, \quad \forall k\geq1,
        $$
        \begin{equation}\label{alphak 5}
             \alpha_{k}(x)\to+\infty \quad \forall x \in  \Omega_{\ell}\setminus\overline{\Omega}_{\ell}^0,
        \end{equation}
        and
        \begin{equation}\label{beta k 5}
    \beta_{k}(x)\to+\infty  \quad \forall x \in  \overline{\Omega}_{n\ell}\setminus\overline{\Omega}_{n\ell}^0.
        \end{equation}
Then, 
\begin{equation}\label{muk-to-Lambda0}
    \mu_{k}  \to \Lambda_0:=\left\{
\begin{array}{ll}
     \lambda_p^{\Omega_0}[\mathfrak{L}_0+A]& \mbox{if } \operatorname{dist}(\Omega_{\ell}^0, \Omega_{n\ell}^0)< \delta,   \\
     \min \left\{\sigma_{1}^{\Omega_{\ell}^0}[-\Delta +a_1],\lambda_p^{\Omega_{n\ell}^0}[L_{n\ell}^{0}+a_2]\right\}& \mbox{if } \operatorname{dist}(\Omega_{\ell}^0, \Omega_{n\ell}^0)\geq  \delta,
\end{array}\right.
\end{equation}
where $\Omega_0:=\Oe^0 \cup \Onl^0$.

Furthermore, 
\begin{enumerate}
    \item[\rm(i)] If $\operatorname{dist}(\Omega_{\ell}^0, \Omega_{n\ell}^0)< \delta$ and $\lambda_p^{\Omega_0}[\mathfrak{L}_0+A]$ admits an eigenfunction $\varphi_0 \in X_{0} $ with $\varphi_0 \gg 0$ and $|\varphi_0|_{L^{2}(\Omega_0)}=1$, then
    $$
    \varphi_k \to \phi \quad \mbox{in } H_{\Gamma_0}^{1}(\Oe)\times L^2(\Onl)
    $$
    where 
      $$\phi:= \left\{ 
        \begin{array}{ll}
            \varphi_0 & \mbox{in }\Omega_{0},  \\
             0& \mbox{in }\Omega\setminus {\Omega}_{0}.
        \end{array}
        \right.$$
    Moreover,  
    $$\varphi_{k} \to \varphi_0 \quad \mbox{in } L_{loc}^{\infty}(\Oe^0)\times C(\overline{\Omega}_{n\ell}^0).$$
    
    \item[\rm(ii)] If $\operatorname{dist}(\Omega_{\ell}^0, \Omega_{n\ell}^0)\geq \delta$ and $\sigma_{1}^{\Omega_{\ell}^0}[-\Delta +a_1]<\lambda_p^{\Omega_{n\ell}^0}[L_{n\ell}^{0}+a_2]$, for $k$ sufficiently large, $\mu_k$ admits an eigenfunction $\varphi_k\in X_\mathfrak{B}$ such that $\varphi_k\gg0$, $|\varphi_k|_2=1$, satisfying
    $$
\varphi_k \to \phi \quad \mbox{in } H_{\Gamma_0}^{1}(\Oe)
$$
where
    $$\phi:= \left\{ 
        \begin{array}{ll}
            \varphi_0 & \mbox{in }\Omega_{\ell}^{0},  \\
             0& \mbox{in }\Omega_{\ell}\setminus \overline{\Omega}_{\ell}^{0},
        \end{array}
        \right.$$
          and $\varphi_0 \in W^{2,p}_0(\Omega_{\ell}^{0})$, for $p>N$, denotes the positive eigenfunction associated with $\sigma_{1}^{\Omega_{\ell}^0}[-\Delta +a_1]$, normalized such that $|\varphi_0|_{L^2(\Omega_{\ell}^{0})}=1$. Furthermore,
        $$
        \varphi_k \to \varphi_0 \quad \mbox{in }L^\infty_{loc}(\Omega_{\ell}^{0}).
        $$
        
\item[\rm(iii)] If $\operatorname{dist}(\Omega_{\ell}^0, \Omega_{n\ell}^0)\geq \delta$, $\lambda_p^{\Omega_{n\ell}^0}[L_{n\ell}^{0}+a_2]<\sigma_{1}^{\Omega_{\ell}^0}[-\Delta +a_1]$, and $\lambda_p^{\Omega_{n\ell}^0}[L_{n\ell}^{0}+a_2]$ admits a positive principal eigenfunction $\varphi_{0} \in C(\overline{\Omega}_{n\ell}^{0})$ with $|\varphi_{0}|_{L^{2}(\Omega_{n\ell}^{0})}=1$, then for $k$ sufficiently large, $\mu_k$ admits an eigenfunction $\varphi_k\in X_\mathfrak{B}$ such that $\varphi_k\gg0$, $|\varphi_k|_2=1$, satisfying
$$
\varphi_k \to  \varphi_0 \quad \mbox{in } C(\overline{\Omega}_{n\ell}^{0}).
$$ 
\end{enumerate} 
\end{proposition}

\begin{proof}
Arguing as in the previous proofs, the sequence $\{\mu_{k}\}$ is non-decreasing, satisfies $\mu_k \leq \Lambda_0$, and the following limit exists:
    \begin{align*}
        \lim_{k\to + \infty}\mu_{k}=\mu^{*} \le \Lambda_{0}.
    \end{align*}
    Suppose by contradiction that $\mu^{*}<\Lambda_{0}$. In particular, we have 
    \begin{equation}
        \label{mu* -cota caso geral}\mu^{*}<\sigma_{1}:=\sigma_{1}^{\Omega_{\ell}^0}[-\Delta +a_1] \quad \mbox{and}\quad \mu^{*}<\lambda_{n\ell}^0:=\lambda_p^{\Omega_{n\ell}^0}[L_{n\ell}^{0}+a_2].
    \end{equation}
    Since
    \begin{align*}
    \mu_{k} \to \mu^*< \lambda_{n\ell}^0\leq \min_{\overline{\Omega}_{n\ell}^0} a_2(x) \quad \mbox{and} \quad\min_{\overline{\Omega}_{n\ell}}(a_{2}(x)+\beta_{k}(x)) \to \min_{\overline{\Omega}_{n\ell}^{0}}a_{2}(x),
    \end{align*}
    it follows that for $k \ge 1$ sufficiently large,
    \begin{align*}
        \mu_{k} < \min_{\overline{\Omega}_{n \ell}}(a_{2}(x)+\beta_{k}(x)).
    \end{align*}
    By Theorem \ref{Th: caso l_p<m} combined with Lemma \ref{lema-varphik-local}, there exists an eigenfunction $\varphi_{k} \in X_{\mathfrak{B}}$ associated with $\mu_{k}$ such that $\varphi_{k} \gg 0$, $|\varphi_{k}|_{2}=1$, satisfying
    \begin{align*}
\varphi_k \rightharpoonup \varphi_\ell \quad \text{in } H_{\Gamma_0}^1(\Oe),\quad \quad 
\varphi_k\to \varphi_{\ell}\quad \text{in } L^2(\Oe),
\end{align*}
and $$\varphi_{\ell}=0 \quad\mbox{a.e. in } \Oe \setminus \overline{\Omega}_{\ell}^{0} \quad\mbox{and } \varphi_{\ell} \in H_{{0}}^{1}(\Omega_{\ell}^{0}).$$

Moreover, since $\mu^{*} < \lambda_{n\ell}^{0} \le \displaystyle\min_{\overline{{\Omega}}_{n\ell}^{0}}(a_2(x))$, we can repeat the arguments from the proof of Proposition \ref{refugio nao local} to deduce that
\begin{equation}  \label{2.51}
    \varphi_{k} \to 0 \quad \text{in } L^{2}(\Onl \setminus \Omega_{n\ell}^{0}) \quad \mbox{and} \quad  \varphi_{k} \to \varphi_{n\ell}  \quad \mbox{in }C(\overline{\Omega}_{n\ell}^{0}),
\end{equation}
where 
$$
\varphi_{n\ell}:= \lim_{k \to +\infty} \frac{-L_{n\ell}\varphi_k-L_{\ell}\varphi_k}{a_2(x)-\mu_k}.
$$
Consequently,
$$
L_{n\ell}\varphi_k = -  2\int_{\Onl} J(x-y) \varphi_k(y) \,dy \to -  2\int_{\Onl^0} J(x-y) \varphi_{n\ell}(y) \,dy=L_{n\ell}^0 \varphi_{n\ell} \quad
\mbox{in } C(\overline{\Omega}_{n\ell}),$$
and
$$
L_{\ell}\varphi_k = -  \int_{\Oe} J(x-y) \varphi_k(y) \,dy \to -  \int_{\Oe^0} J(x-y) \varphi_{\ell}(y) \,dy=L_{\ell}^0 \varphi_{\ell} \quad
\mbox{in } C(\overline{\Omega}_{n\ell}).
$$

Now, multiplying the equation for $\varphi_k|_{\Oe}$ by $\psi \in H_{0}^1(\Omega_{\ell}^{0})$ and integrating by parts over $\Oe^0$, we deduce that
\begin{multline*}
     \int_{\Omega_{\ell}^{0}}\nabla \varphi_{k}\nabla \psi+\int_{\Omega_{\ell}^{0}}a_{1}(x)\varphi_{k}\psi -\int_{\Omega_{\ell}^{0}}\psi \int_{\Onl^{0}}J(x-y)\varphi_{k}(y)dydx  \\-\int_{\Omega_{\ell}^{0}}\psi \int_{\Onl \setminus \Onl^{0}}J(x-y)\varphi_{k}(y)dydx=\mu_{k}\int_{\Omega_{\ell}^{0}}\varphi_{k}\psi.
\end{multline*}
On the other hand, $\varphi_{k}|_{\Onl^0}$ satisfies
\begin{equation*}
    L_{n\ell}\varphi_{k} + L_{\ell}\varphi_{k} + a_2(x)\varphi_{k}  = \mu_{k}\varphi_{k} \quad \text{in } \Onl^0.
\end{equation*}
Letting $k \to +\infty$ and using the previous convergences, it becomes apparent that 
\begin{multline}\label{eq-local0}
    \int_{\Omega_{\ell}^{0}}\nabla \varphi_{\ell}\nabla \psi+\int_{\Omega_{\ell}^{0}}a_{1}(x)\varphi_{\ell}\psi -\int_{\Omega_{\ell}^{0}}\psi \int_{\Onl^{0}}J(x-y)\varphi_{n\ell}(y)dydx \\=\mu^*\int_{\Omega_{\ell}^{0}}\varphi_{\ell}\psi \quad \forall \psi \in H_{0}^1(\Oe^0),
\end{multline}
and
\begin{equation}\label{eq-naolocal0}
    L_{n\ell}^0\varphi_{n\ell} + L_{\ell}^0\varphi_{\ell} + a_2(x)\varphi_{n\ell}= \mu^*\varphi_{n\ell} \quad \text{in } \Onl^0.
\end{equation}

Assume that $\operatorname{dist}(\Omega_{\ell}^0,\Omega_{n\ell}^0)<  \delta$. Since $\varphi_\ell$ satisfies \eqref{eq-local0}, by an elliptic regularity argument, we have $\varphi_{\ell} \in W^{2,p}(\Omega_{\ell}^{0})$ $(p>N)$. Thus, equations \eqref{eq-local0}-\eqref{eq-naolocal0} imply that $\varphi_0:=(\varphi_{\ell},\varphi_{n\ell}) \in W_{0}^{2,p}(\Oe^{0}) \times C(\overline{\Omega}_{n\ell}^{0}) $ is a nonnegative and nontrivial solution of
\[
\begin{cases} 
\mathfrak{L}_{0}u +A(x)u=\mu^*u & \text{in } \Omega_{0},\\ 
 u  = 0& \text{on } \partial \Omega_{\ell}^{0}.

\end{cases}
\]
By the uniqueness of the principal eigenvalue, it follows that $\mu^*=\lambda_{p}^{\Omega_{0}}[\mathfrak{L}_{0}+A]$, which contradicts the hypothesis.

Suppose now that $\operatorname{dist}(\Omega_{\ell}^0,\Omega_{n\ell}^0)\geq  \delta$. It follows that $B_\delta(x)\cap\Onl^{0} = \emptyset$ for all $x \in \Oe^0$ and $B_\delta(x)\cap\Oe^{0} = \emptyset$ for all $x \in \Onl^0$. Thus, hypothesis \ref{J1} implies that
$$ \int_{\Onl^{0}}J(x-y)\varphi_{n\ell}(y)dy =\int_{B_\delta(x)\cap\Onl^{0}}J(x-y)\varphi_{n\ell}(y)dy = 0 \quad \forall x \in \Omega_{\ell}^0$$
and
$$
L_{\ell}^0\varphi_{\ell}(x)= - \int_{\Oe^0} J(x-y) \varphi_{\ell}(y)\, dy =  - \int_{B_\delta(x)\cap\Oe^{0}} J(x-y) \varphi_{\ell}(y)\, dy=0 \quad \forall x \in \Omega_{n\ell}^0.
$$
Therefore, \eqref{eq-local0}-\eqref{eq-naolocal0} are equivalent to
\begin{align*}
    \int_{\Omega_{\ell}^{0}}\nabla \varphi_{\ell}\nabla \psi+\int_{\Omega_{\ell}^{0}}a_{1}(x)\varphi_{\ell}\psi =\mu^*\int_{\Omega_{\ell}^{0}}\varphi_{\ell}\psi \quad \forall \psi \in H_{0}^{1}( \Omega_{\ell}^{0}).
\end{align*}
and
\begin{equation*}    L_{n\ell}^0\varphi_{n\ell} + a_2(x)\varphi_{n\ell}= \mu^*\varphi_{n\ell} \quad \text{in } \Onl^0.
\end{equation*}

Since
\begin{align*}
        \int_{\Oe^0} \varphi_\ell^2+\int_{\Onl^0} \varphi_{n\ell}^2 &= \lim_{k \to +\infty}\left(\int_{\Omega_{\ell}^0}\varphi_{k}^{2}+\int_{\Oe \setminus \Omega_{\ell}^{0}}\varphi_{k}^{2}+\int_{\Omega_{n\ell}^0}\varphi_{k}^{2}+\int_{\Onl \setminus \Omega_{n\ell}^{0}}\varphi_{k}^{2}\right)
        \\&= \lim_{k \to +\infty} |\varphi_k|_2^2  = 1,
\end{align*}
it follows that $\varphi_\ell,\varphi_{n\ell} \geq 0$ and either $\varphi_{\ell} \neq 0$ in $\Oe^0$ or $\varphi_{n\ell} \neq 0$ in $\Onl^0$. By the uniqueness of the principal eigenvalue, the above identities imply that $\mu^*= \sigma_1$ or $\mu^*=\lambda_{n\ell}^0$, which contradicts \eqref{mu* -cota caso geral}. This concludes the proof of \eqref{muk-to-Lambda0}. 

In what follows, we will prove (i). If $\lambda_{p}^{\Omega_{0}}[\mathfrak{L}_0+A]$ admits a positive eigenfunction $\varphi_0 \in X_{0}$ normalized in $L^2(\Omega_0)$, by Proposition \ref{Prop:lambda-av-principal} we have 
$$\lambda_{p}^{\Omega_{0}}[\mathfrak{L}_{0}+A] < \min_{\overline{\Omega}_{n\ell}^{0}}(a_{2}(x)).$$
Since $$\mu_{k} \to \lambda_{p}^{\Omega_{0}}[\mathfrak{L}_0+A] \quad \mbox{and } m_{k}:=\min_{\overline{\Omega}_{n\ell}}(a_{2}(x)+\beta_{k}(x))\to \min_{\overline{\Omega}_{n\ell}^{0}}(a_{2}(x)),$$
for $k \ge 1$ sufficiently large, we have $\mu_k< m_k$. Therefore, by Theorem \ref{Th: caso l_p<m}, there exists an eigenfunction $\varphi_{k} \in X_{\mathfrak{B}}$ associated with $\mu_k$ such that $\varphi_{k} \gg 0$ and $|\varphi_{k}|_{L^{2}(\Omega)}=1$. Define
 $$\phi:= \left\{ 
        \begin{array}{ll}
            \varphi_0 & \mbox{in }\Omega_{0},  \\
             0& \mbox{in }\Omega\setminus {\Omega}_{0}.
        \end{array}
        \right.$$
Repeating the arguments above, we obtain
$$
\varphi_k \to \phi \quad \mbox{in } L^2(\Omega),
$$
$$
\varphi_k \rightharpoonup \phi \quad \text{in } H_{\Gamma_0}^1(\Oe) \quad \mbox{and} \quad \varphi_k \to \varphi_0 \quad \mbox{in } C(\overline{\Omega}_{n\ell}^0).$$
To show that 
$$\varphi_{k} \to \varphi_{\ell}:= \phi|_{\Oe} \quad \mbox{in } H^{1}(\Oe),$$
we proceed exactly as in the proof of Proposition \ref{Refugio em local} to conclude that 
$$\normaind{\nabla \varphi_{k}}{L^{2}(\Oe)} \to \normaind{\nabla \varphi_\ell}{L^{2}(\Oe)}.$$
In particular, since $\Oe^{0} \subset \Oe$, Sobolev embeddings imply that 
$$\varphi_{k} \to \varphi_0 \quad \mbox{in } L^{2^*}(\Oe^{0}).$$ 
Finally, by an elliptic regularity argument (see Proposition \ref{Refugio em local}), we obtain 
$$\varphi_{k} \to \varphi_0 \quad \mbox{in } L^{\infty}_{loc}(\Oe^{0}).$$
Consequently,
 $$\varphi_{k} \to \varphi_0 \quad \mbox{in } L_{loc}^{\infty}(\Oe^0)\times C(\overline{\Omega}_{n\ell}^0).$$
To prove case (ii), it suffices to observe that $\sigma_{1}^{\Omega_{\ell}^0}[-\Delta +a_1]<\lambda_p^{\Omega_{n\ell}^0}[L_{n\ell}^{0}+a_2] \le \min_{\overline{\Omega}_{n\ell}^{0}}(a_{2}(x))$. Hence,
$$\mu_{k} \to \sigma_{1}^{\Oe^{0}}[-\Delta+a_{1}] \quad \mbox{and } m_{k}:=\min_{\overline{\Omega}_{n\ell}}(a_{2}(x)+\beta_{k}(x))\to \min_{\overline{\Omega}_{n\ell}^{0}}(a_{2}(x)).$$
Thus, we can repeat the arguments from Proposition \ref{Refugio em local} to conclude the result.

Finally, let us prove item (iii). Suppose that $\lambda_{p}^{\Omega_{n\ell}^{0}}[L_{n\ell}^{0}+a_2]$ admits a positive principal eigenfunction $\varphi_{0}$ normalized in $L^{2}(\Onl^{0})$. By Lemma \ref{propriedades-NL-autovalor}, we have $\lambda_{p}^{\Omega_{n\ell}^{0}}[L_{n\ell}^{0}+a_2]<\min_{\overline{\Omega}_{n\ell}^0}(a_{2}(x))$. 

Since $\operatorname{dist}(\Omega_{\ell}^0, \Omega_{n\ell}^0)\geq \delta$ and $\lambda_p^{\Omega_{n\ell}^0}[L_{n\ell}^{0}+a_2]<\sigma_{1}^{\Omega_{\ell}^0}[-\Delta +a_1]$, we deduce that 
 $$\mu_{k} \to \lambda_{p}^{\Omega_{n\ell}^{0}}[L_{n\ell}^{0}+a_2] \quad \mbox{and } m_{k}:=\min_{\overline{\Omega}_{n\ell}}(a_{2}(x)+\beta_{k}(x))\to \min_{\overline{\Omega}_{n\ell}^{0}}(a_{2}(x)).$$
Consequently, $\mu_k< m_k$ for $k \ge 1$ sufficiently large. Therefore, by using the arguments from Proposition \ref{refugio nao local}, we conclude the proof.
\end{proof}

\begin{remark}
    We emphasize that in Proposition \ref{refugio local e nlocal}, the hypothesis $\overline{\Omega}_{\ell}^0 \subset \Oe$ can be replaced by $\Oe^0 \subset \Oe$ with $\Gamma \subset \partial \Oe^0$ or $\Gamma
    \cap \partial \Oe^0 = \emptyset$. The same applies to Proposition \ref{Refugio em local}.
\end{remark}

Note that Theorem \ref{T2} follows directly from Propositions \ref{Refugio em local}, \ref{refugio nao local}, and \ref{refugio local e nlocal}.

\section{Behavior with respect to \texorpdfstring{$\delta$}{delta}}\label{sec5}

The main goal of this section is to establish the behavior of $\lambda_{p}^\Omega[\mathfrak{L}+A]$ with respect to the parameter $\delta$. To this end, recall that hypothesis \ref{J1} states that $\operatorname{supp} J = B_\delta(0)$. Thus, in this section, it is convenient to use the notation 
$$J_\delta=J \quad \text{and} \quad \mathfrak{L}_\delta:=\mathfrak{L}.$$ 
For the sake of simplicity, we will also write
$$\lambda(\delta) := \lambda_{p}^{\Omega}[\mathfrak{L}_{\delta}+A].$$

The following result establishes that if $J_\delta$ is monotonic with respect to $\delta$, then $\lambda(\delta)$ will also be monotonic.

\begin{proposition}\label{monotonia com respeito a delta}
    Assume that the map $\delta \in (0,+\infty) \mapsto J_\delta \in C(\mathbb{R}^N)$ is increasing, in the sense that $J_{\delta_1}(z) \le J_{\delta_2}(z)$ for all $z \in \mathbb{R}^N$ whenever $0 < \delta_1 < \delta_2$. Then the map $\delta \in (0,\infty) \mapsto \lambda(\delta)$ is decreasing. That is, if $0 < \delta_{1} < \delta_{2}$, then 
    $$\lambda(\delta_{2}) \le \lambda(\delta_{1}).$$
\end{proposition}

\begin{proof}
    Define 
\begin{equation*}
    \mathcal{A}_{\delta_1} := \Big\{ \lambda \in \mathbb{R} : \exists \, \phi \in X, \, \phi(x) > 0, \, x \in \Omega, \text{ supersolution of } (\mathfrak{L}_{\delta_1} + A- \lambda, \mathfrak{B}, \Omega) \Big\}
\end{equation*}
and 
\begin{equation*}
    \mathcal{A}_{\delta_{2}} := \Big\{ \lambda \in \mathbb{R} : \exists \, \phi \in X, \, \phi(x) > 0, \, x \in \Omega, \text{ supersolution of } (\mathfrak{L}_{\delta_2} + A - \lambda, \mathfrak{B}, \Omega) \Big\}.
\end{equation*}
By Proposition \ref{autovalor com supersoluçao}, it suffices to show that 
$$\mathcal{A}_{\delta_{2}} \subset  \mathcal{A}_{\delta_1}.$$
Let $\lambda \in \mathcal{A}_{\delta_{2}}$. Then, there exists $\phi \in X$ with $\phi(x) > 0$ for $x \in \Omega$ such that 
\begin{align}
    \begin{cases}
\mathfrak{L}_{\delta_2}\phi+A(x)\phi -\lambda \phi  \ge 0& \text{in } \Omega,\\ 
\mathfrak{B} \phi \ge   0& \text{on } \partial \Omega_{\ell}.
\end{cases} \label{supersolução em delta2}
\end{align}
On the other hand, combining the fact that $J_{\delta}$ is increasing with $B_{\delta_1}(x) \subset B_{\delta_2}(x)$, we find that 
\begin{align*}
-\int_{\Omega_{n\ell} \cap B_{\delta_1}(x)} J_{\delta_1}(x-y)\phi(y) \, dy 
&\ge -\int_{\Omega_{n\ell} \cap B_{\delta_1}(x)} J_{\delta_2}(x-y)\phi(y) \, dy \\
&\ge -\int_{\Omega_{n\ell} \cap B_{\delta_2}(x)} J_{\delta_2}(x-y)\phi(y) \, dy.
\end{align*}
Similarly, 
\begin{align*}
-\int_{\Omega_{\ell} \cap B_{\delta_1}(x)} J_{\delta_1}(x-y)\phi(y) \, dy 
&\ge -\int_{\Omega_{\ell} \cap B_{\delta_1}(x)} J_{\delta_2}(x-y)\phi(y) \, dy \\
&\ge -\int_{\Omega_{\ell} \cap B_{\delta_2}(x)} J_{\delta_2}(x-y)\phi(y) \, dy.
\end{align*}
Substituting the above inequalities into \eqref{supersolução em delta2}, it becomes apparent that
\begin{align*}
    \begin{cases}
\mathfrak{L}_{\delta_1}\phi+A(x)\phi -\lambda \phi  \ge \mathfrak{L}_{\delta_2}\phi+A(x)\phi -\lambda \phi \ge  0& \text{in } \Omega,\\ 
\mathfrak{B} \phi \ge   0& \text{on } \partial \Omega_{\ell}.
\end{cases}
\end{align*}
Thus, we find that $\lambda \in \mathcal{A}_{\delta_{1}}$. Consequently, from the properties of the supremum, we infer that 
$$\lambda(\delta_{2}) = \lambda_{p}^{\Omega}[\mathfrak{L}_{\delta_{2}}+A] \le \lambda_{p}^{\Omega}[\mathfrak{L}_{\delta_{1}}+A] = \lambda(\delta_{1}).$$
This concludes the proof.
\end{proof}

Now we establish the asymptotic behavior of $\lambda(\delta)$ as $\delta \to 0^+$.

\begin{proof}[Proof of Theorem \ref{T3}]
To simplify the notation, let us denote $\sigma_1 := \sigma_{1}^{\Oe}[-\Delta+a_{1}; \mathfrak{B}]$ and $m := \min_{\overline{\Omega}_{n\ell}}a_2(x).$ 

By Proposition \ref{estimativas} and inequality \eqref{lambdap<m}, it follows that 
$$\lambda(\delta) \le \min \{\sigma_1, m\} \quad \forall \delta > 0.$$
If \eqref{delta-to-0} did not hold, there would exist a sequence $\delta_k > 0$, with $k \geq 1$, and a constant $c < \min\{\sigma_1, m\}$ such that
$$
\delta_k \to 0  \quad \text{and} \quad \lambda(\delta_k) < c,
$$
for all $k \geq 1$. In particular,
    \begin{align}
        \lambda(\delta_{k}) < \min\{\sigma_1, m\} \le m. \label{autovalor delta k menor que a2}
    \end{align}
    
    Thanks to the strict inequality in \eqref{autovalor delta k menor que a2}, it follows from Theorem \ref{Th: caso l_p<m} that there exists an eigenfunction $\varphi_{k} \in X_{\mathfrak{B}}$ associated with $\lambda(\delta_k)$ such that $\varphi_{k} \gg 0$ and $|\varphi_{k}|_2 = 1$. Multiplying the equation of $\varphi_k|_{\Oe}$ by $\psi \in H_{\Gamma_0}^1(\Oe)$ and integrating by parts over $\Oe$ yields
        \begin{align}
       \int_{\Oe}\nabla \varphi_{k} \cdot \nabla\psi+\int_{\Oe} a_{1}(x) \varphi_{k} \psi-\int_{\Oe} \psi \int_{\Onl \cap B_{\delta_{k}}(x)}J_{\delta_{k}}(x-y)\varphi_{k}(y) \, dydx=\lambda(\delta_{k}) \int_{\Oe} \varphi_{k} \psi.\label{lambdak e varphk}
    \end{align}
    Since $\delta_{k} \to 0^{+}$, for sufficiently large $k \geq 1$ we have $\delta_{k} < {\delta}'$. Choosing $\psi=\varphi_{k}|_{\Oe}$, we can use Hölder's inequality combined with \eqref{J-limitada} to conclude that there exists a constant $\tilde{C}>0$ such that
    $$\lambda(\delta_{k}) \ge -\tilde{C}-\normaind{a_{1}}{\infty}.$$
     Consequently $\{\lambda(\delta_k)\}$ is bounded and, passing to a subsequence still denote by $k$, we obtain
    $$
    \lambda(\delta_k)\to \lambda^* \le c.
    $$
    Notice that \eqref{lambdak e varphk} also implies that $\{\varphi_k\}$ is bounded in $H_{\Gamma_0}^1(\Oe)$. Thus, for a subsequence still labeled by $k$, we have
    $$\varphi_{k} \rightharpoonup \varphi_{\ell} \quad \text{in } H_{\Gamma_0}^1(\Oe) \quad \text{and} \quad  \varphi_{k} \to \varphi_{\ell} \text{ in } L^{2}(\Oe).$$
    In particular, since $a_2(x) - \lambda(\delta_k) > 0$ for all $x \in \Onl$ and $k \geq 1$, we deduce that $\varphi_{k}|_{\Onl}$ satisfies 
    \begin{align}
        \varphi_{k}=\frac{-L_{n\ell}\varphi_{k}- L_{\ell}\varphi_{k}}{a_{2}(x)-\lambda(\delta_{k})} \quad \text{in }\Onl.\label{igualdade para varphik com deltak}
    \end{align}
     We claim that 
     $$L_{n\ell}\varphi_{k}(x)=-2\int_{\Omega_{n\ell} \cap B_{\delta_k}(x)} J_{\delta_k}(x-y)\varphi_{k}(y) \, dy \to 0 \quad \text{uniformly in }  \overline{\Omega}_{n\ell}.$$
      Combining the fact that $|\varphi_{k}|_2 = 1$ with Hölder's inequality and hypothesis \eqref{J-limitada}, it follows that
      \begin{align*}
          |L_{n\ell}\varphi_k(x)|&\le 2\int_{\Omega_{n\ell} \cap B_{\delta_k}(x)} |J_{\delta_k}(x-y)\varphi_{k}(y)| \, dy \nonumber \\ 
          & \le 2\left(\int_{\Onl \cap B_{\delta_k}(x)} |J_{\delta_k}(x-y)|^2 \,dy\right)^\frac{1}{2}\\
          &\le2 C |\Onl \cap B_{\delta_k}(x)|^{1/2}\\
          &\le2 C |B_{\delta_k}(0)|^{1/2} \quad \forall x \in \overline{\Omega}_{n\ell}. \label{desigualdade com a caracteristica}
      \end{align*}
      Consequently,
      $$\|L_{n\ell}\varphi_k\|_{C(\overline{\Omega}_{n\ell})} \to 0.$$
      Similarly, we also obtain
      $$\|L_{\ell}\varphi_k\|_{C(\overline{\Omega}_{n\ell})} \to 0.$$
      Since 
      $$
      a_2(x) - \lambda(\delta_k) \to a_2(x) - \lambda^* > 0 \quad \forall x \in \overline{\Omega}_{n\ell},      $$
      letting $k \to \infty$ in \eqref{igualdade para varphik com deltak}, we find that 
      $$\varphi_{k} \to 0 \quad\text{in } C(\overline{{\Omega}}_{n\ell}). $$
      Subsequently, using the convergence results established above, we can let $k \to \infty$ in \eqref{lambdak e varphk} to deduce that 
      $$\int_{\Oe}\nabla \varphi_\ell \cdot \nabla\psi+\int_{\Oe} a_{1}(x) \varphi_{\ell} \psi=\lambda^{*} \int_{\Oe} \varphi_{\ell} \psi \quad \forall \psi \in H_{\Gamma_0}^1(\Oe).$$
      Since 
      $$  1=\lim_{k \to \infty}|\varphi_{k}|_{2}^{2}=\lim_{k \to \infty}\left(\int_{\Omega_\ell}\varphi_{k}^{2}+\int_{\Onl}\varphi_{k}^{2}\right)=|\varphi_{\ell}|^2_{L^2(\Omega_\ell)},$$
      we have that $\varphi_{\ell} \ge 0$ and $\varphi_{\ell} \not\equiv 0$. 
      Consequently, we deduce that $\lambda^* = \sigma_{1}$, contradicting our initial assumption that $\lambda^* < \min\{\sigma_1, m \} \le \sigma_1$. This concludes the proof.

\end{proof}

\begin{remark}
We point out that the conclusion of Theorem \ref{T3} remains valid if we assume that $\operatorname{dist}(\Oe,\Onl)=d_0>0$. Indeed, in this case, problem \eqref{P} decouples for $0<\delta\leq d_0$. Consequently, by the definition of $\lambda_p^\Omega[\mathfrak{L}+A]$ (see \eqref{Def lambdap}), we obtain
\begin{equation}
    \label{lambda(delta)-caso2}
    \lambda(\delta) = \min \left\{\sigma_{1}^{\Oe}[-\Delta+a_{1};\mathfrak{B}], \lambda_p^{\Onl}[L_{\delta}+a_2]\right\} \quad \forall \delta \in (0,d_0],
\end{equation}
where $L_\delta = L_{n\ell}$. Now, using \eqref{lambdap-NL<m} and taking $u\equiv 1$ in the definition of $\lambda_p^{\Onl}[L_{\delta}+a_2]$ (see \eqref{autovalor nonlocal}), we have
\begin{align*}
    \min_{\overline{\Omega}_{nl}} a_2(x)&\ge \lambda_p^{\Omega_{nl}}[L_{\delta}+a_2]\\
    &\ge \inf_{x \in \Omega_{nl}} \left( a_2(x) - 2\int_{\Omega_{nl} \cap B_\delta(x)} J_\delta(x-y)\,dy \right) \\
    &\ge \min_{\overline{\Omega}_{nl}} a_2(x) - 2 C |B_\delta(0)|, \quad \forall \delta \in (0, \delta')
\end{align*}
    where $\delta'$ is given by  hypothesis \eqref{J-limitada}.
Therefore,
\[
\lambda_p^{\Omega_{nl}}[L_{\delta}+a_2] \to \min_{\overline{\Omega}_{nl}} a_2(x) \quad \text{as } \delta \to 0^+.
\]
Thus, letting $\delta \to 0^+$ in \eqref{lambda(delta)-caso2}, we recover \eqref{delta-to-0}.    
\end{remark}

\subsection*{Funding}

\noindent W. C. has been  supported by CNPq Proc. 310664/2023-3.  J. P. was financed in part by CAPES Proc. 88881.125193/2025-01 and 88887.328317/2026-00. A. S. has been supported under Grants PID2023-149509NB-I00 and  CEX2024-001517-M - Apoyo a Unidades de Excelencia María de Maeztu  funded by MICIU/AEI. A. S. also has been supported by University of Sevilla under grants SOL2024-31708 and SOL2024-31596. 

\bibliographystyle{acm} 
\bibliography{myref}
\end{document}